\RequirePackage{scrlfile} 
\PreventPackageFromLoading{natbib}

\documentclass[]{elsarticle}
\date{February 1, 2020}
\usepackage{lineno,hyperref}
\modulolinenumbers[5]
\journal{Computers \& Mathematics with Applications}
\usepackage{a4wide}
\usepackage[utf8]{inputenc}
\usepackage{amsmath}
\usepackage{amsfonts}
\usepackage{amssymb}
\usepackage{amsthm}
\usepackage{mathtools}
\usepackage{bm} %bolt math font
\usepackage{stmaryrd} %for double brackets
\usepackage[USenglish]{babel}
\usepackage{caption}
\usepackage{cleveref}
\usepackage{epstopdf}
\usepackage{subcaption}% for subfigures
\usepackage{eso-pic}%title page figure placement
\usepackage{siunitx} %for scientific numbers \num
\usepackage{cite}

\definecolor{pscol}{rgb}{0,0,0.7}
\definecolor{ipcol}{rgb}{0.7,0,0}
\newcommand{\ps}[1]{#1}%{\color{pscol}{#1}}}
\usepackage{tikz}
\tikzset{>=stealth} % set arrows as stealth fighter jets
\usepackage{tkz-euclide}
\usetikzlibrary{positioning,arrows,calc,decorations.markings,math,arrows.meta}
\usepackage{multirow}
\usepackage[normalem]{ulem}
\useunder{\uline}{\ul}{}

\theoremstyle{definition}
\newtheorem{definition}{Definition}[section]

\newtheorem{remark}[definition]{Remark}
\newtheorem*{remark*}{Remark}

\newtheorem*{acknowledgment*}{Acknowledgment}

\theoremstyle{plain}

\numberwithin{equation}{section}

\newcommand{\R}{\mathbb{R}}

\newcommand\restr[2]{{% we make the whole thing an ordinary symbol
\left.\kern-\nulldelimiterspace % automatically resize the bar with \right
#1 % the function
\right|_{#2} % this is the delimiter
}}

\newcommand{ \al } { \{\!\!\{ }
\newcommand{ \ar } { \}\!\!\} }
\newcommand{ \jl } { [\![ }
\newcommand{ \jr } { ]\!] }
\newcommand{\dt}[1]{\frac{\partial {#1}}{\partial t}}
\newcommand{\dtt}[1]{\frac{\partial^2 {#1}}{\partial t^2}}

\begin{document}

\begin{frontmatter}

\title{Tent pitching and Trefftz-DG method for the acoustic wave equation}

%% include affiliations in footnotes:
\author[univie]{Ilaria Perugia}
\ead{ilaria.perugia@univie.ac.at}

\author[tuwien]{Joachim Schöberl}
\ead{joachim.schoeberl@tuwien.ac.at}

\author[univie]{\texorpdfstring{Paul Stocker\corref{corrauth}}{Paul Stocker}}
\cortext[corrauth]{Corresponding author}
\ead{paul.stocker@univie.ac.at}

\author[tuwien]{Christoph Wintersteiger}
\ead{christoph.wintersteiger@tuwien.ac.at}

\address[univie]{Faculty of Mathematics, University of Vienna, Austria}
\address[tuwien]{Institute for Analysis and Scientific Computing, TUWien, Austria}

\begin{abstract}
We present a space-time Trefftz discontinuous Galerkin method for approximating the acoustic wave equation semi-explicitly on tent pitched meshes. 
DG Trefftz methods use discontinuous test and trial functions, which solve the wave equation locally. 
Tent pitched meshes allow us to solve the equation elementwise, allowing locally optimal advances in time.
The method is implemented in NGSolve, solving the space-time elements in parallel, whenever possible.
Insights into the implementational details are given, including the case of propagation in heterogeneous media. 
%We are able to recover optimal convergence rates in terms of mesh size and in terms of polynomial degree of the Trefftz spaces. 
\end{abstract}

\begin{keyword}
wave equation, discontinuous Galerkin method, Trefftz method, tent pitched mesh 
\MSC[2010] 65M60\sep 41A10\sep 35L05
\end{keyword}

\end{frontmatter}

\sloppy
\section{Introduction}
%This is a study on the combination of the space-time Trefftz discontinuous Galerkin (DG) method \cite{Moiola2017} with tent pitched meshes \cite{tp1,tp2}.\par
Standard finite element methods approximate the solution of a given partial differential equation (PDE) by piecewise polynomial functions. 
A classic approach to the discretization of time dependent PDEs is to use finite element methods to discretize space and then use time stepping schemes to advance in time. 
We consider a different approach based on finite element approximation simultaneously in space and time. 
\ps{In this paper, we present original numerical results for the
  space-time Trefftz discontinuous Galerkin (DG) method studied in \cite{Moiola2017},
confirming the known, and conjectured, properties of the method numerically.
This is the first implementation combining a Trefftz-DG method with
tent pitched meshes, highlighting the importance of both techniques
in order to produce an efficient method.
We point out that we can
%are able to
solve the problem in $n+1$ dimensions, $n\in\{1,2,3\}$,
without the need for $n+1$ %four
dimensional elements, since the Trefftz method only %just
requires the computation of %surface
integrals at interelement boundaries.
Furthermore, we introduce a way of recovering the solution of the second order system, and address in detail the issue of inhomogeneous materials.
}
The space-time approach requires us to mesh the full space-time domain. 
Furthermore, the use of approximation spaces based on piecewise "total degree" polynomials in both space and time, leads to a higher number of degrees of freedom.
On the upside, $hp$-refinement is made possible in space-time, allowing for straightforward higher order approximation. 
Additionally, the space-time domain mesh is not forced to be a product mesh, as it is for time stepping schemes. 
Instead, we are allowed to use unstructured meshes. 
This gives us the possibility to devise suitable mesh design strategies in order to circumvent the CFL-condition, which usually limits the global time-step size by the size of the smallest spatial element in explicit time stepping schemes. 
This will allow to advance in time elementwise and in parallel.
% both are innovative techniques in the discretization of wave propagation problems.
%These can be suitably designed in order to take into account possible causality constraints. 
\medskip\par
Space-time finite element methods for linear wave propagation go as far back as \cite{Hughes}, and have been used with DG methods e.g. in \cite{Pete-2005,FalkRichter,LILIENTHAL,Drfler}.
DG methods are based on discontinuous piecewise polynomial functions, and so-called numerical fluxes which impose continuity constraints at mesh inter-element boundaries.
\par
On the one hand, ideas on combining DG methods with so called tent pitched meshes can be found in e.g. \cite{Richter,Pete-2005,FalkRichter,Lowrie,10.1007/978-3-642-59721-3_48}.
A way of constructing these meshes can, for instance, be found in \cite{tp1,tp2,MTP,WinterMTP}.
Tent pitching techniques generate a space-time mesh, which complies with the causality properties of the hyperbolic PDE.  
The resulting mesh consists of tent shaped objects, each advancing
locally optimally in time, with the PDE being explicitly solvable in
each of them.
Though the tent pitching strategy pairs well with DG methods, also other methods are applicable in combination with tent pitched meshes.
In \cite{MTP,WinterMTP,MTPMaxwell}, schemes for the semi-discretization of different hyperbolic equations on tents are presented, which map tents to a domain where space and time are separable. 
Similar to the Trefftz-DG method, these schemes are able to solve in 3+1 dimensions, without building four dimensional elements.
%This allows to have locally optimal time advances in the mesh.
% obtaining error estimates of order $O(h^{p+\frac12})$, for polynomials of degree $p$.
% On the other hand,
Friedrichs theory is used in \cite{gopala} in order to derive a
conforming methods,
and to prove its convergence properties.
We point out that tent pitching is not the only way to deal with the time step restriction of locally refined meshes.
A stabilization for a conforming space-time finite element method on Cartesian (in time) meshes is presented in \cite{Steinbach2019}.
Classical time-stepping schemes can still be applied successfully by splitting the domain into a coarse-mesh and a fine-mesh region, then explicit time stepping in the coarse-mesh region is combined with local implicit or explicit time stepping in the fine-mesh region.
A fully explicit scheme can be found in \cite{grote1,grote2}.
\par
On the other hand, Trefftz methods, originating from \cite{trefftz}, incorporate properties of the PDE into the test and trial spaces.
This is done by choosing them as (local) solutions of the targeted differential equation.
The use of Trefftz spaces allows us to reduce the number of degrees of freedom, as compared to the total degree polynomial spaces, however keeping the same accuracy.
Work on Trefftz-DG methods for different wave propagation problems includes \cite{Kocher,Banjai,barucq,EggerKretzMaxwell,KretzMoiolaMaxwell,Moiola2017,Egger2015TransparentBC}.
In \cite{Banjai}, a Trefftz-DG method in space-time for the second order wave equation is presented, proving $h$-convergence in 1, 2 and 3 space dimensions, as well as $hp$-convergence, along with exponential convergence for analytic solutions in 1 space dimension.
In \cite{barucq}, Trefftz-DG is applied to the coupled elasto-acoustic system, and well-posedness of the problem, as well as error estimates in mesh-dependent norms, are shown. 
Both \cite{Banjai, barucq} are formulated for meshes with tensor product structure in time.
In \cite{Moiola2017} a Trefftz-DG method for the acoustic wave equation in first order formulation and in arbitrary space dimension is presented.
The formulation works for tensor product (in time) meshes, as well as for tent pitched meshes. 
Well-posedness and optimal $h$-convergence are proven, \ps{but no
numerical results are shown.} 
\par
In this paper, we focus on the combination of Trefftz-DG formulation
and tent pitched meshes, and on its efficient implementation.
As Trefftz-DG formulations only contain interelement terms, they pose a natural choice to evolve the solution from the bottom to the top of tent elements. 
\ps{As mentioned, previous implementations on Trefftz-DG methods were
  based on tensor product structure in time \cite{Banjai,barucq,KretzMoiolaMaxwell}.}
\par This work proceeds as follows.
First, we introduce the Trefftz-DG method in \Cref{sec:themethod}, starting by stating the model problem, defining the Trefftz spaces and finishing the section by formulating the method, as it was introduced in \cite{Moiola2017}.
We continue in \Cref{sec:Trefftzdisc} by reviewing different strategies of discretizing the Trefftz spaces. 
In \Cref{sec:evoltent}, we discuss some numerical details on how to evolve the solution elementwise on a tent pitched mesh, and in 
\Cref{sec:recovery}, we show a way to recover the second order solution from the first order formulation. 
Finally, we present numerical results, which were obtained by implementation of the method in NGSolve \cite{joachim,ngsolve}, in \Cref{sec:numerics}.

\section{The Trefftz-DG method for the acoustic wave equation}\label{sec:themethod}
\subsection{The acoustic wave equation}
Let $ Q=\Omega\times(0,T)$ be a space-time domain in $\R^{n+1},\ n\in\{1,2,3\}$, where $\Omega\subset\R^n$ is a Lipschitz bounded domain with outward unit normal vector $\bm n^{\bm x}_\Omega\in\R^n$.
\ps{The boundary of $\Omega$ is divided into two parts $\Gamma_D$
  and $\Gamma_N$, corresponding to Dirichlet and Neumann boundary
  conditions, respectively, such that they have disjoint
  interiors, the union of their closures gives the boundary of
    $\Omega$, and one of them can be empty.}
We consider the acoustic wave equation in first order formulation, given by
\begin{align}\label{eq:firstordwaveeq}
    \begin{cases}
        \nabla \cdot \bm\sigma + c^{-2} \frac{\partial v}{\partial t} = 0  & \qquad \text{ in } Q,\\
        \nabla v + \frac{\partial \bm\sigma}{\partial t} = \bm 0 & \qquad \text{ in } Q, \\
        v(\cdot,0) = v_0,\ \bm\sigma(\cdot,0) = \bm\sigma_0 & \qquad \text{ on } \Omega, \\
        v=g_D & \qquad \text{ on } \Gamma_D \times [0,T], \\
        \bm n^x_\Omega \cdot \bm\sigma = g_N & \qquad \text{ on } \Gamma_N \times [0,T], \\
    \end{cases}
\end{align}
where we assume that the wavespeed $c>0$ is piecewise constant on $\Omega$. 
\par If the initial condition $\bm\sigma_0$ it the gradient of a scalar field $U_0$, i.e. $\bm\sigma_0=-\nabla U_0$, then  the first order system is equivalent to the second order system obtained by setting $v=\dt{U}$ and $\bm\sigma=-\nabla U$:
\begin{align}\label{eq:secondordwaveeq}
    \begin{cases}
        -\Delta U+ c^{-2} \dtt U = 0 & \qquad \text{ in } Q, \\
        \dt{U}(\cdot,0)=v_0,\ U(\cdot,0)=U_0 & \qquad \text{ on } \Omega, \\
        \dt{U}=g_D & \qquad \text{ on } \Gamma_D \times [0,T], \\
        -\bm n^x_\Omega \cdot \nabla U = g_N & \qquad \text{ on } \Gamma_N \times [0,T]. 
    \end{cases}
\end{align}
The Laplacian $\Delta$, gradient $\nabla$ and divergence $\nabla\cdot$ are considered with respect to the space variable $\bm x$ only.

\subsection{Space-time meshes}
The mesh $\mathcal{T}_h(Q)$ of the space-time domain $Q$ is assumed to consist of non-overlapping Lipschitz polytopes, where $h=\max_{K\in\mathcal{T}_h(Q)}h_K$, with $h_K$ being the anisotropic diameter defined in \eqref{eq:adiam}.
For each mesh face $F=\partial K_1\cap \partial K_2$, for $K_1,K_2\in\mathcal{T}_h(Q)$, we assume that it either lies below the characteristic speed $1/c$, or is vertical (parallel to the time axis). 
In more rigorous terms:
Let $(\bm n_F^x,n_F^t)$ be the normal vector to $F$ with $n_F^t\geq 0$, then either 
\begin{align*}
    & c|\bm n_F^x|<n^t_F \quad \text{ and we call the face \textit{space-like}, or} \\
    & n^t_F=0 \quad \text{ and we call the face \textit{time-like}.}
\end{align*}
Notice, however, that no CFL-condition or any other time step size restriction is imposed on the time-like faces.
\par
A mesh with space-like faces only is called a tent pitched mesh.
\ps{In the numerical experiments presented below, we generate tent
  pitched meshes by using the algorithm
  presented in~\cite{MTP} (see also~\cite{WinterMTP}).}
The mesh is buit by progressively advancing in time, stacking tent-shaped objects on top or each other, each of them union of $(n+1)$-simplices. %, $n$ being the space-dimension. 
The main idea is that the tent height is chosen such that the differential equation is explicitly solvable in each tent. 
Therefore, the local maximal time advance at a spatial point has to respect the causality constraint, which corresponds to a \textit{local} CFL-condition. 
This allows us to advance the solution tent by tent, not necessarily having to solve a global system. 
For independent tents, i.e. tents that are not on top of each other,
the computations can be done in parallel.
\ps{We remark that, in the numerical tests, we observe no stability issues
 with tents pitched very close to the limit of the causality
 condition.}

%The main benefit of such a mesh is that the solution on each element can be obtained explicitly, if the initial condition on the bottom (lower in time) faces of the element.  This requirement can be seen as a local CFL-condition. 

%%%%%%%%%%%%%%%%%%%%%%%%%%%%%%%%%%%%%%%%%%%%%%%%%%%%%%%%%%%%%%%%%%%%%%%%%%%%%%%%%%%%%%%%%
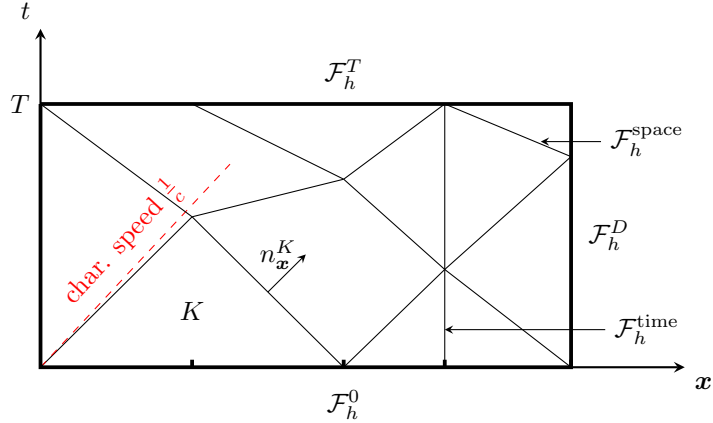
\begin{figure}[!ht]
        \centering
    \resizebox{0.6\linewidth}{!}{
    \begin{tikzpicture}
        %coord. system
        \draw[thick,->] (0,0) -- (8.5,0) node[anchor=north west] {$\bm x$};
        \draw[thick,->] (0,0) -- (0,4.5) node[anchor=south east] {$t$};
        %gray tent
        % \fill[fill=gray!30] (7,2.8) -- (5+1/3,3.5) -- (7,3.5); %1
        % 0
        \draw[] (0,0) -- (2,2);			 % 1
        \draw[] (4,0) -- (2,2);
        \draw[] (4,0) -- (5+1/3,1.3);		 % 2
        \draw[] (7,0) -- (5+1/3,1.3); 		%3
        \draw[red,dashed] (0,0) -- (2+.5,2.2+.5) node [midway, above, sloped]{\ \ char. speed $\frac 1c$};			 % 1
        % 1
        \draw[] (2,2) -- (0,3.5);		 % 1
        \draw[] (2,2) -- (4,2.5);		 % 2
        \draw[] (5+1/3,1.3) -- (4,2.5);
        \draw[] (5+1/3,1.3) -- (7,2.8);	 %3
        %2
        \draw[] (4,2.5) -- (5+1/3,3.5) ; %2
        \draw[] (7,2.8) -- (5+1/3,3.5) ;
        \draw[] (4,2.5) -- (2,3.5) ; %1
        %draw thick boundary
        \draw[line width=0.5mm] (0,3.5)node[anchor= east] {$T$} -- (0,0) -- (7,0) -- (7,3.5) -- (0,3.5);
        %adv front
        %\draw[thick] (0,3.5) -- (7,3.5) ; %1
        %normal
        \draw[->] (3,1) -- (3.5,1.5)   node[anchor= east] {$n_{\bm x}^K$};
        %labels
        \node at (2,0.75) {$K$};
        \node at (4,-0.5) {$\mathcal{F}_h^0$};
        \node at (4,3.9) {$\mathcal{F}_h^T$};
        \node at (7.5,1.75) {$\mathcal{F}_h^D$};
        \node at (8,3) (space-like) {$\mathcal{F}_h^{\text{space}}$};
        \draw[->] (space-like) -- (6.6,3);
        \draw[] (5+1/3,0) -- (5+1/3,3.5);
        \node at (8,0.5) (time-like) {$\mathcal{F}_h^{\text{time}}$};
        \draw[->] (time-like) -- (5+1/3,0.5);
        %tics
        \draw[line width=0.5mm] (2,0) -- (2,.1);
        \draw[line width=0.5mm] (4,0) -- (4,.1);
        \draw[line width=0.5mm] (5+1/3,0) -- (5+1/3,.1);
        \draw[line width=0.5mm] (7,0) -- (7,.1);
    \end{tikzpicture}
    }
    \caption{\ps{A tent pitched mesh with faces below the characteristic speed and a time-like faces, which are contained inside two tents}.} 
    \label{fig:mesh}
\end{figure}
%%%%%%%%%%%%%%%%%%%%%%%%%%%%%%%%%%%%%%%%%%%%%%%%%%%%%%%%%%%%%%%%%%%%%%%%%%%%%%%%%%%%%%%%
\par The set of all faces of our mesh,
\begin{align*}
    \mathcal{F}_h & :=\bigcup_{K\in\mathcal{T}_h}\partial K,
\end{align*}
can be separated into two sets of internal faces
\begin{align*}
    \mathcal{F}_h^{\text{space}} & := \{F\in\mathcal{F}_h : F \text{ is internal space-like face}\},\\
    \mathcal{F}_h^{\text{time}} & := \{F\in\mathcal{F}_h : F \text{ is internal time-like face}\},
\end{align*}
and sets of faces on the boundary of $Q$, split up according to their types of initial and boundary conditions
\begin{align*}
    \mathcal{F}_h^0 & :=\Omega\times\{0\}, \quad \mathcal{F}_h^T := \Omega\times\{T\}, \\
    \mathcal{F}_h^D &:=\Gamma_D\times[0,T], \quad \mathcal{F}_h^N:=\Gamma_N\times[0,T].
\end{align*}
This classification of the faces is represented in \Cref{fig:mesh}.
\subsection{Trefftz spaces}
By definition, Trefftz functions in the kernel of the considered differential operator. 
For the first order wave equation, we define the local and global Trefftz space as
\begin{align*}
    \bm T(K):=\big\{ 
    (w,\bm\tau)\in L^2(K)^{1+n} \text{ s.t. } \restr{\bm\tau}{\partial K}\in L^2(\partial K)^n,\ \dt w,\nabla\cdot\bm\tau\in L^2(K),\\
    \dt{\bm\tau},\nabla w\in L^2(K)^{1+n},\  
    \nabla \cdot \bm\tau + c^{-2} \frac{\partial w}{\partial t} = 0,\   
    \nabla w + \frac{\partial \bm\tau}{\partial t} = \bm 0\
    \big\}
\end{align*}
\begin{align*}
    \bm T(\mathcal T_h):=\Big\{
    (w,\bm\tau)\in L^2(Q)^{1+n},\ \text { s.t. } (\restr{w}{K},\restr{\bm\tau}{K})\in\bm T(K)\ \forall K\in\mathcal T_h
    \Big\},
\end{align*}
respectively.
Note that, by assuming that the solution is in $\bm T(\mathcal T_h)$, we require additional smoothness on the solution, as in general we only have that $\dt v + \nabla\cdot\bm\sigma\in L^2(K)$, for all $K\in\mathcal T_h$.
\par We derive the Trefftz-DG method for any choice of discrete test and trial space with a Trefftz property, which we denote by $\bm V_p(\mathcal T_h)$.
A possible choice for a polynomial $\bm V_p(\mathcal T_h)$ is given in \Cref{sec:Trefftzdisc} below.

\subsection{The Trefftz-DG method}
Following \cite{Moiola2017}, we derive the Trefftz-DG method for the IBVP in \eqref{eq:firstordwaveeq}.
The method is derived from a local weak formulation, obtained by multiplying the two equations in \eqref{eq:firstordwaveeq} by test and trial functions $w$ and $\bm\tau$, respectively, and integrating by parts on each element $K$ of the mesh $\mathcal{T}_h(Q)$. 
Then, adding the two equations gives
\begin{align}\label{eq:weakform}
    \begin{split}
        & -\int_K v\left(\nabla\cdot\bm\tau+c^{-2}\dt w\right) + \bm\sigma\cdot\left(\dt{\bm\tau}+\nabla w\right)\ dV \\
        & +\int_{\partial K} v \left( \bm\tau \cdot\bm n^x_K + c^{-2} w n^t_K\right) +   {\bm\sigma}\cdot \left(w \cdot\bm n^x_K + \bm\tau n^t_K  \right)\  dS = 0.
    \end{split}
\end{align}
By choosing Trefftz test functions $(w,\bm\tau)\in\bm{V}_p(K)$, the volume integrals over $K$ vanishes. 
We are left with:
\begin{align}\label{eq:Trefftzweakform}
    \begin{split}
        & \int_{\partial K} \hat{v}_{hp} \left( \bm\tau \cdot\bm n^x_K + c^{-2} w n^t_K\right) +   \hat{\bm\sigma}_{hp}\cdot \left(w \cdot\bm n^x_K + \bm\tau n^t_K  \right)\  dS = 0.
    \end{split}
\end{align}
Typical for DG methods, the continuity of the numeric solution on inter-element boundaries is enforced within the bilinear form of the method. 
To this end, the trace of the solution $(v,\bm\sigma)$ in the boundary integral has been replaced by the numeric fluxes $( \hat{v}_{hp},\hat{\bm\sigma}_{hp})$, which we define below.

\medskip To do so, we need to introduce some standard DG notation. 
For a face $F=\partial K^+\cap \partial K^-$ shared by two elements $K^+,K^-\in\mathcal{T}_h(Q)$, we define the average $\al \cdot \ar$, jumps in space $\jl\cdot\jr_N$, and jumps in time $\jl\cdot\jr_t$, for scalar- and vector-valued functions as follows
\begin{align*}
    \al w\ar &:=\frac12 (\restr{w}{K^+}+\restr{w}{K^-}),&  
    \al \bm\tau\ar &:=\frac12 (\restr{\bm\tau}{K^+}+\restr{\bm\tau}{K^-})  \\
    \jl w\jr_t &:= \restr{w}{K^+} n_{K^+}^t + \restr{w}{K^-} n_{K^-}^t, &
    \jl \bm\tau \jr_t &:= \restr{\bm\tau}{K^+} n_{K^+}^t + \restr{\bm\tau}{K^-} n_{K^-}^t,\\
    \jl w\jr_N &:= \restr{w}{K^+}\bm n_{K^+}^x + \restr{w}{K^-}\bm n_{K^-}^x, &
    \jl \bm\tau \jr_N &:= \restr{\bm\tau}{K^+}\cdot \bm n_{K^+}^x + \restr{\bm\tau}{K^-}\cdot \bm n_{K^-}^x,
\end{align*}
where $n_K=(\bm n_K^x,n_K^t)$ is the unit outer normal vector at $\partial K$ split into its space and time components. 
%For scalar functions, the dot product in the definition of the jump in space is replaced by a scalar product. 

\medskip 
%Recall that the mesh $\mathcal{T}_h(Q)$ of the space-time domain can be either a product mesh in time or a tent pitched mesh.
%On the product mesh we need to fix the information loss across the time-like faces, using the averages and jumps defined above. 
Across time-like faces, the information is passed by using \ps{centered fluxes with jump penalization}, whereas, across space-like faces, the information is passed upward in time, resembling an up-wind scheme. 
More precisely, the fluxes on the inter-element faces are chosen as
$$
\hat{v}_{hp} = \begin{cases}
    v^-_{hp} \\
    v_0 \\
    \al v_{hp}\ar + \beta \jl\bm\sigma_{hp}\jr_N \\
    v_{hp} \\
    g_D \\
    v_{hp}+\beta(\bm\sigma_{hp}\cdot\bm n^x_\Omega-g_N)
\end{cases}
\hat{\bm\sigma}_{hp}=\begin{cases}
    \bm\sigma^-_{hp} & \quad \text{on } \mathcal{F}_h^{\text{space}}, \\
    \bm\sigma_0 & \quad \text{on } \mathcal{F}_h^{0},\\
    \al \bm\sigma_{hp}\ar + \alpha \jl v_{hp}\jr_N & \quad \text{on } \mathcal{F}_h^{\text{time}},\\
    \bm\sigma_{hp} & \quad \text{on } \mathcal{F}_h^{T},\\
    \bm\sigma_{hp} - \alpha(v_{hp}-g_D) \bm n^x_\Omega & \quad \text{on } \mathcal{F}_h^{D},\\
    g_N\bm n^x_\Omega & \quad \text{on } \mathcal{F}^N_h,
\end{cases}
$$
where $\alpha$ and $\beta$ are penalty parameters, which will be chosen constant (notice that they are needed on time-like and Dirichlet faces only).
By $w^+$ and $w^-$ we denote the trace of the function $w$ on space-like faces from the adjacent element at higher and lower times, respectively.
% $\alpha=\beta=0\rightarrow$ [Egger \& al., 2014]\ \\[.4cm]
% $\alpha\beta\geq \frac14\rightarrow$ [Monk, Richter, 2005]
%
% \begin{figure}
% 	\includegraphics[width=.7\linewidth]{img/DGspace}
% 	\label{fig:dgspace}

Finally, we plug the definition of the fluxes into \eqref{eq:Trefftzweakform} and sum over all elements $K\in\mathcal{T}_h(Q)$. Then the Trefftz-DG method for the wave equation reads:
\begin{align}\begin{split}\label{eq:TDG}
    & \text{find } (v_{hp},\bm\sigma_{hp})\in\bm{V}_p(\mathcal{T}_h) \quad\text{s.t.} \\
    & \mathcal{A}(v_{hp},\bm\sigma_{hp};w,\bm\tau) = \ell(w,\bm\tau) \qquad \forall (w,\bm\tau)\in\bm{V}_p(\mathcal T_h),
\end{split}\end{align}
with
\begin{align*}  
    \mathcal{A}(v_{hp},\bm\sigma_{hp};w,\bm\tau ) &:=  \int_{\mathcal{F}_h^{\text{space}}} \left( c^{-2} v_{hp}^- \jl w\jr_t +\bm\sigma^-_{hp}\cdot \jl\bm\tau\jr_t + v^-_{hp} \jl\bm\tau\jr_N+\bm\sigma^-_{hp} \cdot \jl w\jr_N \right)\ dS \\
    & + \int_{\mathcal{F}_h^{\text{time}}} \left( \al v_{hp}\ar \jl \bm\tau\jr_N + \al\bm\sigma_{hp}\ar\cdot\jl w\jr_N + \alpha \jl v_{hp}\jr_N\cdot\jl w \jr_N + \beta \jl\bm\sigma_{hp}\jr_N\jl\bm\tau\jr_N \right)\ dS \\
    & + \int_{\mathcal{F}^T_h} c^{-2}v_{hp} w + \bm\sigma_{hp}\cdot\bm\tau\ dS + \int_{\mathcal{F}^D_h} (\bm\sigma\cdot\bm n^x_\Omega+\alpha v_{hp} ) w \ dS 
 \\ & + \int_{\mathcal F_h^N} v_{hp}(\bm\tau\cdot\bm n^x_\Omega)+\beta(\bm\sigma\cdot\bm n_\Omega^x)(\bm\tau\cdot\bm n^x_\Omega)\ dS
\end{align*}
and
\begin{align*}
    \ell(w,\bm \tau) &:= \int_{\mathcal{F}_h^0} c^{-2} v_0 w + \bm\sigma\cdot\bm\tau \ dS + \int_{\mathcal{F}_h^D} g_D(\alpha w - \bm\tau\cdot\bm n^x_\Omega)\ dS  + \int_{\mathcal F_h^N} g_N(\beta\bm\tau\cdot\bm n^x_\Omega - w)\ dS.
\end{align*}

\medskip\par 
%On the product mesh the time-slabs can be solved sequentually, on each of which the method is implicit.
On a tent pitched mesh, as the one in \Cref{fig:mesh}, the method is semi-explicit, meaning that the solution on each tent only depends on the tents below, allowing to solve each tent explicitly, and tents independent from each other in parallel; details are given in \Cref{sec:evoltentconst} below.
The situation where also vertical faces are present, is needed, for instance, in the case of piecewise constant wavespeed, is discussed in \Cref{sec:evoltentpiece} below.
%Of course one could also solve the whole tent slab at once, implicitly.
Note that the method only includes integrals over element boundaries, thus only quadrature on $n$ dimensional simplices is needed.

\section{Choice of discrete Trefftz spaces}\label{sec:Trefftzdisc}
So far, we have not specified what discretization of the Trefftz space $\bm V_p(K)\subset \bm T(K)$ to use.
We introduce the straightforward choice, given by all polynomials in space-time that fulfill the first order wave equation.
For an element $K\subset\mathbb{R}^{n+1}$ in the mesh $\mathcal{T}_h(Q)$, we define the local polynomial Trefftz space as
\begin{align*}
    \mathbb{T}^p(K)&:=\mathbb P^{p}(\mathbb R^{n+1})^{n+1}\cap\bm T(K),
\end{align*}
where we denote by $\mathbb{P}^p(K)$ the space of polynomials on $K$ of degree $\leq p$.
In general, it is possible to choose different polynomial degrees $p$ in different elements. 
Here, we choose a uniform $p$, as this is consistent with the numerical examples below.
The global Trefftz-DG space on the whole mesh is then given by $\mathbb{T}^p(\mathcal{T}_h):= \prod_{K\in \mathcal{T}_h} \mathbb{T}^p(K).$
%The main property of the Discontinuous Galerkin method is that the global basis functions are piecewise continuous (thus non-conforming), as seen in the choice of the global approximation space.
The dimension of the elemental Trefftz space is given by  
\begin{equation*}
    \dim\mathbb{T}^p(K)=(n+1)\begin{pmatrix}p+n\\ n\end{pmatrix}=\mathcal{O}_{p\rightarrow\infty} (p^n),
\end{equation*}
where we recall that
$\left(\begin{smallmatrix}
        a\\ b
\end{smallmatrix}\right)=\frac{a!}{b!(a-b)!}$ for $b\leq a\in\mathbb{N}_0$.
Notice that, for the total degree polynomial space, one has $\dim(\mathbb{P}^p(\mathbb R^{n+1})^{n+1})=\mathcal{O}_{p\rightarrow\infty} (p^{n+1})$.
\medskip\par Let us now assume that the first order problem is derived from a second order problem.
Then it is natural to derive the vector valued Trefftz space for the first order problem from a scalar Trefftz space for the second order problem. 
We now detail this approach as it is the one we use for the numerical results presented in \Cref{sec:numerics}.
Let us start by defining the polynomial Trefftz space for the second order problem: 
\begin{equation*}
    \mathbb{U}^p(K):=\left\{U\in \mathbb{P}^p(K): -\Delta U+\frac{1}{c^2}\frac{\partial^2}{\partial t^2} U=0 \right\}.
\end{equation*}
We are able to construct a basis for this space using the recursion formula introduced in \cite[Remark 13]{Moiola2017}. 
We recall it here, for completeness.
We need some multi-index notation: for $\bm\alpha\in\mathbb{N}^n_0$ we denote $|\bm\alpha|=\alpha_1+\dots+\alpha_n$ and $\bm x^{\bm\alpha}=x_1^{\alpha_1}\dots x_n^{\alpha_n}$. 
Furthermore, let $\bm e_m:=(0,\dots,0,1,0,\dots,0)\in\mathbb{N}^n_0$ with 1 in the $m$-th entry. 
Consider a space-time polynomial 
\begin{equation*}
    U(\bm x,t)=\sum_{\substack{\bm\alpha\in\mathbb N^n_0,\ k\in\mathbb N_0,\\ |\bm\alpha|+k\leq p}} a_{k,\bm\alpha}\bm x^{\bm\alpha}t^k.
\end{equation*}
We want to compute the coefficients $a_{k,\bm\alpha}$ such that the polynomial is Trefftz.
This is done by inserting the polynomial into the second order wave equation and collecting terms of equal power to find that 
\begin{equation}\label{eq:Trefftz2ndorderrecursion} 
    a_{k,\bm\alpha} = \frac{c^2}{k(k-1)} \sum_{m=1}^n(\alpha_m+1)(\alpha_m+2)a_{k-2,\bm\alpha+2\bm e_m}
\end{equation}
has to hold for the polynomial to be Trefftz.
To start the recursion, we need to choose polynomial bases (in the space variables only) for $k=0$ and $k=1$, respectively. 
More precisely, we start by choosing polynomial basis functions $ \{\tilde{b}_1,\dots,\tilde{b}_{\left(
            \begin{smallmatrix}
                p+n\\ n
            \end{smallmatrix}\right)}\} $ for the space $\mathbb P^p(\mathbb R^n)$ and $\{\hat{b}_1,\dots,\hat b_{\left(
            \begin{smallmatrix}
                p-1+n\\ n
\end{smallmatrix}\right)}\}$ for $\mathbb P^{p-1}(\mathbb R^n)$.
Then we can introduce a basis for $\mathbb U^p(K)$ such that either $U(\cdot,0)=\tilde b_\ell$ and $\dt U(\cdot,0)=0$, or $U(\cdot,0)=0$ and $\dt U(\cdot,0)=\hat b_\ell$ for some $\ell$.
Hence, we can construct the basis for $\mathbb U^p(K)$ out of two sets of polynomial basis functions of $\mathbb P^p(K)$ and $\mathbb P^{p-1}(K)$. This lets us determine the dimension as 
\begin{equation*} 
    \dim\mathbb U^p(K)=\begin{pmatrix}p+n\\ n\end{pmatrix} + \begin{pmatrix}p-1+n\\ n\end{pmatrix} = \frac{2p+n}{p}\begin{pmatrix}p-1+n\\ n\end{pmatrix} .
\end{equation*}
\par Then, a Trefftz space for the first order system can be derived from 
\begin{equation*}
    \mathbb U^{p+1}(K)=\text{span}\{b_j,\ j\in \mathcal I\} \text{ by setting } \mathbb W^{p}(K)=\text{span}\{(\dt {b_j},-\nabla b_j),\ j\in\mathcal I\}.
\end{equation*}
We have that 
\begin{equation*}
    \dim \mathbb W^p(K)= \dim\mathbb U^{p+1}(K) -1=\frac{2p+n+2}{p+1}
    \begin{pmatrix} n+p\\ n
    \end{pmatrix}
\end{equation*}
and furthermore $\mathbb W^p(K)\subset \mathbb T^p(K)$.
A recursion formula, similar to \eqref{eq:Trefftz2ndorderrecursion}, can also be derived for $\mathbb T^p(K)$, however the numerical results in \Cref{sec:numerics} are centered around $\mathbb W^p(K)$. 

\begin{remark}\label{rem:scalepol}
    It is sufficient to compute the coefficients only once for $c=1$ and then fix the wavenumber by a coordinate transform. 
    Furthermore, for numerical stability, it is convenient to shift the basis functions to the center of the element and scale them by its anisotropic diameter, which is defined by 
    \begin{equation}\label{eq:adiam}
        h_K:=\sup_{(\bm x,t),(\bm y,s)\in K} (|\bm x - \bm y|^2+c^2|t-s|^2)^{1/2}
    \end{equation}
    for a mesh element $K$. 
    For reference coordinates $(\hat{\bm x},\hat t)$, the coordinate transform given by 
    $$(\bm x,t)=(h_K\hat{\bm x},h_Kc^{-1}\hat t) $$
    transforms the Trefftz basis $\hat U(\hat{\bm x},\hat t)$ of wavespeed 1 to Trefftz basis functions $\hat U(x,t)$ of arbitrary wavespeed $c$. 
    In the case of Trefftz functions for the first order system $(\hat v,\hat{\bm\sigma})$, we need to choose $$v(\bm x ,t)=c\hat v(\hat{\bm x},\hat t),\quad \bm\sigma(\bm x,t)=\hat{\bm\sigma}(\hat{\bm x},\hat t).$$
\end{remark}

\section{Evolution within a tent}\label{sec:evoltent}
The tent pitched mesh allows us to solve local tents explicitly.
This is due to the fact that the slope of the mesh faces is below the characteristic speed $1/c$, thus the local solution on a tent can be computed once the solution on its inflow boundary is known.
In \Cref{sec:evoltentconst}, we discuss how to evolve the solution within a tent with constant wavespeed inside the tent itself. 
The case where the wavespeed changes within a tent is considered in \Cref{sec:evoltentpiece}.
Notice that, in the constant wavespeed case, tents coincide with mesh elements, while in the latter case tents on the interface contain more than one mesh element.

\subsection{Constant wavespeed} \label{sec:evoltentconst}
%\begin{sloppypar}
    Let us denote the bottom and top faces of the tent by ${T^{\text{bot}}_h\subset( \mathcal F_h^{\text{space}}\cup\mathcal F_h^0)}$ and ${T^{\text{top}}_h\subset (\mathcal F_h^{\text{space}}\cup\mathcal F_h^T)}$, respectively. 
    Furthermore, tent faces on the boundary are denoted by $T^D_h\subset \mathcal F^D_h$ for Dirichlet and $T^N_h \subset \mathcal F^N_h$ for Neumann boundaries.
    \par Since the solution is explicit on each tent, we only need to solve a local system  of size ${\dim(\bm V^p(K))\times\dim(\bm V^p(K))}$. The system is derived from \eqref{eq:TDG} and is given by the following equation 
%\end{sloppypar}
\begin{equation}\label{eq:tentsys}
    \begin{aligned}
    & \int_{T^{\text{top}}_h} c^{-2}  v_{hp} w n_K^t + \bm\sigma_{hp}\cdot\bm\tau n_K^t + v_{hp}\bm\tau\cdot \bm n_K^x + \bm\sigma_{hp}\cdot(w\bm n_K^x)\ dS
 \\ &\quad +  \int_{T^D_h} (\bm\sigma_{hp}\cdot\bm n^x_\Omega+\alpha v_{hp} ) w \ dS + \int_{T^N_h} v_{hp}(\bm\tau\cdot\bm n^x_\Omega)+\beta(\bm\sigma\cdot\bm n_\Omega^x)(\bm\tau\cdot\bm n^x_\Omega)\ dS
 \\
    & = -\int_{T^{\text{bot}}_h}  c^{-2} v_{\text{bot}} w n_K^t + \bm\sigma_{\text{bot}}\cdot\bm\tau n_K^t + v_{\text{bot}}\bm\tau\cdot \bm n_K^x + \bm\sigma_{\text{bot}}\cdot\bm n_K^x w\ dS
  \\&  \quad + \int_{T^D_h} g_D(\alpha w - \bm\tau\cdot\bm n^x_\Omega)\ dS + \int_{T^N_h} g_N(\beta\bm\tau\cdot\bm n^x_\Omega - w)\ dS,
    \end{aligned}
\end{equation}
where, in the case $T_{h}^{\text{bot}}\subset\mathcal{F}_h^0$, $(v_{\text{bot}},\bm\sigma_{\text{bot}})=(v_0,\bm\sigma_0)$, and in the case $T_h^{\text{bot}}\subset\mathcal{F}_h^{space}$, $(v_{\text{bot}},\bm\sigma_{\text{bot}})$ on a given face is the previously computed solution in the tent sharing that face in lower time. %$=(v^-,\bm\sigma^-)$.  
%,which is known, since we are solving sequentially. 

%%%%%%%%%%%%%%%%%%%%%%%%%%%%%%%%%%%%%%%%%%%%%%%%%%%%%%%%%%%%%%%%%%%%%%%%%%%%%%%%%%%%%%%%%
\begin{figure}[ht]\centering
    \resizebox{.50\linewidth}{!}{
        \begin{tikzpicture}
            %coord. system
            \draw[thick,->] (0,0) -- (4.5,0) node[anchor=north west] {$\bm x$};
            \draw[thick,->] (0,0) -- (0,3.5) node[anchor=south east] {$t$};
            \draw[] (2,1) coordinate (A) -- (3,2) coordinate (B);
            \draw[] (2,1) -- (1,2) coordinate (C);
            \draw[] (2,3) coordinate (top) -- (3,2) ; 
            \draw[] (2,3) -- (1,2) ;
            %tics
            \draw[line width=0.5mm] (1,0) -- (1,.1);
            \draw[line width=0.5mm] (2,0) -- (2,.1);
            \draw[line width=0.5mm] (3,0) -- (3,.1);
            \draw[line width=0.3mm] (1.2,0) -- (1.2,.1);
            \draw[line width=0.3mm] (1.8,0) -- (1.8,.1);
            \draw[line width=0.3mm] (1.5,0) -- (1.5,.1);
            \draw[line width=0.3mm] (2.2,0) -- (2.2,.1);
            \draw[line width=0.3mm] (2.8,0) -- (2.8,.1);
            \draw[line width=0.3mm] (2.5,0) -- (2.5,.1);

            \foreach \x in {2.2,2.5,2.8}
            {
                \draw (\x,0) coordinate (D);
                \draw[dashed] (D) -- (intersection of A--B and D--[shift={(0,100)}]D) coordinate (E); 
                \fill [] (E) circle (1pt); 
                \node[fill=black,scale=0.4] at (intersection of B--top and D--[shift={(0,100)}]D) [rectangle,draw] () {};
            }
            \foreach \x in {1.2,1.5,1.8}
            {
                \draw (\x,0) coordinate (D);
                \draw[dashed] (D) -- (intersection of A--C and D--[shift={(0,100)}]D) coordinate (E); 
                \node[fill=black,scale=0.4] at (intersection of C--top and D--[shift={(0,100)}]D) [rectangle,draw] () {};
                \fill [] (E) circle (1pt); 
            }
            \draw[] (2,1) -- (2,3);
            \node at (3.5,2.5) (timelike) {\footnotesize${T^{\text{time}}_h}$};
            \draw[->] (timelike) -- (2,2.2);
            \node at (3,3) (tin) {\footnotesize${T^{\text{top}}_h}$};
            \draw[->] (tin) -- (2.3,2.7);
            \node at (3.9,2) (tout) {\footnotesize${T^{\text{bot}}_h}$};
            \draw[->] (tout) -- (2.6,1.6);
            \draw[] (2,2) node[anchor=west] {\scriptsize${c=c_1}$};
            \draw[] (2,2) node[anchor=east] {\scriptsize${c=c_2}$};
        \end{tikzpicture}
    }
    \caption{The spatial integration points are mapped to the faces of the tent. The solution is determined using the known input on the bottom integration points (dots), and is evaluated on the top integration points (squares).}
    \label{fig:numevoltent}
\end{figure}
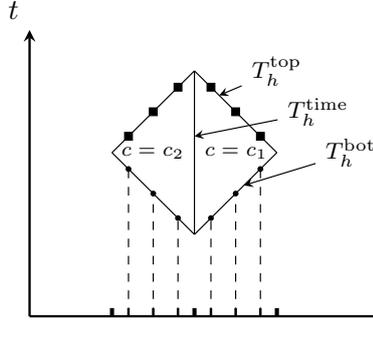
%%%%%%%%%%%%%%%%%%%%%%%%%%%%%%%%%%%%%%%%%%%%%%%%%%%%%%%%%%%%%%%%%%%%%%%%%%%%%%%%%%%%%%%%
\par For the numerical integration, we only need an integration rule for $n$-simplices, in order to integrate over the boundary of the tent. 
We can define an integration rule on the spatial mesh once, which we can then map to the faces of the tent.
This idea is visualized in \Cref{fig:numevoltent}.
After solving on the tent, we need to evaluate $(v_{hp},\bm\sigma_{hp})$ in the integration points on $T^{\text{top}}_h$, and store these values for the next tent. 
On each spatial integration point, we only need to store the most recent results, leading to a total storage of: ${(\text{total number of integration points}) \cdot(n+1)}$. 

\subsection{Piecewise constant wavespeed}\label{sec:evoltentpiece}
Recall that we assume that the wavespeed is constant in time and piecewise constant in space. %and therefore only jumps across inter-element boundaries are possible. 
In this case, we always consider initial spatial meshes that are aligned with the discontinuities of the wavespeed.
To treat such a jump within a space-time tent, we need to incorporate the jump terms from our DG formulation \eqref{eq:TDG}.
This involves integrating on the time-like inter-element boundary contained inside the tent, denoted by $T^{\text{time}}_h\subset \mathcal F^{\text{time}}_h$. 
%We gain the following term on the right hand side of the local system on the tent \eqref{eq:tentsys},
According to \eqref{eq:TDG}, one has to add to the left-hand side of \eqref{eq:tentsys} the term
\begin{equation*}
    \begin{aligned}
        %& \int_{T^{\text{top}}_h} c^{-2}  v_{hp} w n_K^t + \bm\sigma_{hp}\cdot\bm\tau n_K^t + v_{hp}\bm\tau\cdot \bm n_K^x + \bm\sigma_{hp}\cdot(w\bm n_K^x)\ dS \\ &\quad + 
        \int_{T^{\text{time}}_h} \left( \al v_{hp}\ar \jl \bm\tau\jr_N + \al\bm\sigma_{hp}\ar\cdot\jl w\jr_N + \alpha \jl v_{hp}\jr_N\cdot\jl w \jr_N + \beta \jl\bm\sigma_{hp}\jr_N\jl\bm\tau\jr_N \right)\ dS.
        %\\ &\quad +  \int_{T^D_h} (\bm\sigma_{hp}\cdot\bm n^x_\Omega+\alpha v_{hp} ) w \ dS + \int_{T^N_h} (v_{hp}(\bm\tau\cdot\bm n^x_\Omega)+\beta(\bm\sigma\cdot\bm n_\Omega^x)(\bm\tau\cdot\bm n^x_\Omega)\ dS
        %\\
        %& = -\int_{T^{\text{bot}}_h}  c^{-2} v_{\text{bot}} w n_K^t + \bm\sigma_{\text{bot}}\cdot\bm\tau n_K^t + v_{\text{bot}}\bm\tau\cdot \bm n_K^x + \bm\sigma_{\text{bot}}\cdot\bm n_K^x w\ dS
        %\\&  \quad + \int_{T^D_h} g_D(\alpha w - \bm\tau\cdot\bm n^x_\Omega)\ dS + \int_{T^N_h} g_N(\beta\bm\tau\cdot\bm n^x_\Omega - w)\ dS
    \end{aligned}
\end{equation*}
Since the tent now includes two mesh elements, the system matrix is now of size ${2\dim(\bm V^p(K))\times 2\dim(\bm V^p(K))}$.
The extension to interfaces between more than two materials follows. 

\section{Recovery of the solution of the second order equation}\label{sec:recovery}
In the case where the problem comes from a second order formulation we can substitute $v=\dt U$ and $\bm\sigma=-\nabla U$ to write the method in terms of test and trial functions from $\mathbb U_{p+1}(K)$. 
Then the method \eqref{eq:TDG} reads:
\begin{align}
    \begin{split}\label{eq:TDG2}
        & \text{find } U_{hp}\in\mathbb{U}_{p+1}(\mathcal{T}_h) \quad\text{s.t.} \\
        & \hat{\mathcal{A}}(U_{hp};V) = \hat{\ell}(V) \qquad \forall V\in\mathbb{U}_{p+1}(\mathcal T_h),
    \end{split}
\end{align}
with
\begin{equation*} \begin{aligned}
& \hat{\mathcal{A}}(U_{hp}; V ) := \mathcal A(\dt{U_{hp}},-\nabla U_{hp};\dt V,-\nabla V) \text{ and } 
%\\ &\quad \int_{\mathcal{F}_h^{\text{space}}} \left( c^{-2} \dt{ U_{hp}^-} \jl \dt V\jr_t +\nabla U^-_{hp}\cdot \jl\nabla V\jr_t - \dt{U^-_{hp}} \jl\nabla V\jr_N-\nabla U^-_{hp} \cdot \jl \dt V\jr_N \right)\ dS \\
%&\quad + \int_{\mathcal{F}_h^{\text{time}}} \left(- \al \dt{ U_{hp}}\ar \jl \nabla V\jr_N - \al\nabla U_{hp}\ar\cdot\jl \dt V\jr_N + \alpha \jl \dt{ U_{hp}} \jr_N\cdot\jl \dt V \jr_N + \beta \jl\nabla U_{hp}\jr_N\jl\nabla V\jr_N \right)\ dS \\
%&\quad + \int_{\mathcal{F}^T_h} c^{-2}\dt{ U_{hp}} \dt V + \nabla U_{hp}\cdot\nabla V\ dS + \int_{\mathcal{F}^D_h} (\alpha\dt{U_{hp}} -\nabla U_{hp}\cdot\bm n _\Omega ) \dt V \ dS \\[.2cm]
& \hat\ell(V) := \ell(\dt V,-\nabla V).
%\int_{\mathcal{F}_h^0} c^{-2} v_0 \dt V + \nabla U_0\cdot\nabla V \ dS + \int_{\mathcal{F}_h^D} g\ \left(\nabla V\cdot\bm n_\Omega+\alpha \dt V\right) \ dS
\end{aligned} \end{equation*}

\ps{The constant basis function does not contribute to %this
  formulation~\eqref{eq:TDG2}, as only derivatives of the
    unknown $U_{hp}$ are present.}
Thus, this formulation produces the same results as \eqref{eq:TDG} with $\bm V_p(\mathcal T_h) = \mathbb W_p(\mathcal T_h)$.
In order to fix the constants and recover the solution to the second order wave equation, we modify the original formulation by adding the additional terms 
$$\int_{\mathcal{F}_h^{\text{space}}} -\jl U_{hp}\jr_tV^+\ dS +\int_{\mathcal{F}_h^{0}} U_{hp}V\ dS$$
to the bilinear form $\hat{\mathcal{A}}(U_{hp}; V ) $, and 
$$\int_{\mathcal{F}_h^{0}} U_0V\ dS$$
to the right hand side $\hat\ell(V)$, where and $U_0(x)=U(x,0)$. Note that these terms preserve the consistency of the formulation.
\par Therefore, when evolving the solution inside a single tent, we need to add $\int_{T^{\text{top}}_h}U_{hp}V$ and $\int_{T^{\text{bot}}_h}U_{\text{bot}}V$ to the left- and right-hand side, respectively, of the formulation discussed in \Cref{sec:evoltent}.

\section{Numerical results}\label{sec:numerics}
\begin{figure*}[ht]
        \centering
    \begin{subfigure}[t]{0.3\textwidth}
        \centering
        \includegraphics[width=\textwidth]{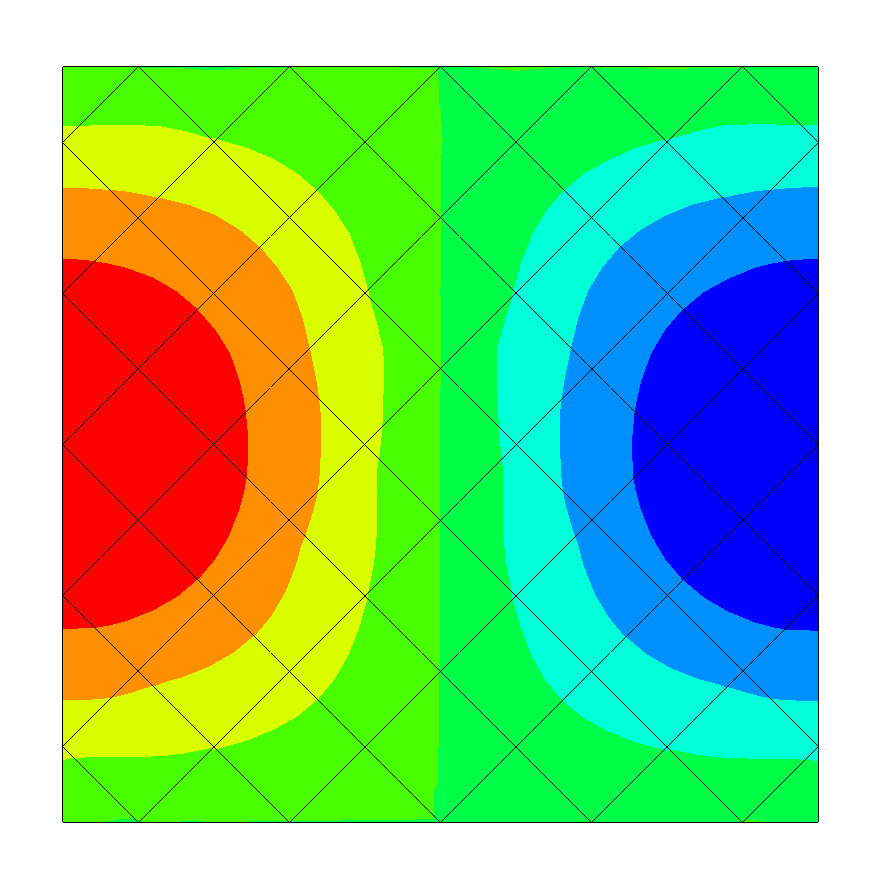}
        %\caption{Conditioning}
    \end{subfigure}
    \caption{Approximation of the standing wave on a 1+1 dimensional space-time tent pitched mesh.}
    \label{fig:tentsol}
\end{figure*}

In this section we present numerical test results in one, two, and three spatial dimensions.
The Trefftz-DG method was implemented in NGSolve \cite{joachim,ngsolve}.  
%on Cartesian meshes (in time) in \Cref{sec:compbasis,sec:comppol}, and on tent pitched meshes in \Cref{sec:tp23,sec:energy,sec:lshape,sec:mat}.
If not otherwise stated, we use the following settings for the numerical examples.
%We compare our numerical solution $(u_{hp},\bm\sigma_{hp})\in \mathbb W^p$ to the following exact analytic solution of the second order wave equation, the standing wave,
We consider the problem \eqref{eq:firstordwaveeq} with initial and Dirichlet boundary conditions such that the analytical solution is $(v,\bm\sigma)=(\dt U,-\nabla U)$, where $U$ is the standing wave
\begin{equation} 
    U(\bm x,t)=\cos(\pi x_1)\cos(\pi x_2)\cos(\pi x_3)\sin(\pi t c \sqrt{n})/(\sqrt{n} \pi), 
\end{equation}
given here in 3+1 dimensions, and set the wavespeed $c=1$. An example is plotted in 1+1 dimensions in \Cref{fig:tentsol}.
%And setting $v = \dt{U}$ and $\bm\sigma=-\nabla U$.
%We choose $v_0=v(\cdot,0)$ and $\bm{\sigma}_0= -\bm\sigma(\cdot,0)$, and use Dirichlet boundary conditions to enforce the solution. 
%We are using the formulation \eqref{eq:TDG2}, employing the Trefftz-space for the second order formulation.
%To recover the second order solution we use the correction terms introduced in the previous section.
The penalty parameters are chosen as $\alpha=\beta=0.5$. 
We measure the error 
\begin{equation*}
        e(v,\bm \sigma;v_{hp},\bm\sigma_{hp})=\left( c^{-2}\| v(\cdot,T)-v_{hp}(\cdot,T)\|^2_{L^2(\Omega)}+ \|\bm\sigma(\cdot,T)-\bm\sigma_{hp}(\cdot,T)\|^2_{L^2(\Omega)}\right)^{\frac12},
\end{equation*}
at final time $T$, which we choose at $T=1$.
\ps{The tent pitched meshes are produced by the algorithm presented in
  \cite{MTP}.
In this algorithm, the height of
 the tents is limited by the slope of the edges, and not by the slope
 of the faces.}
All timings were performed on a server with two Intel(R) Xeon(R) CPU E5-2687W v4, with 12 cores each.

\subsection{Approximation properties of Trefftz spaces}\label{sec:comppol}
Independently of their combination with tent pitching, Trefftz polynomial spaces possess good approximation properties for wave solutions.
In \cite[Section 6]{Moiola2017} , $h$-version approximation estimates for wave solutions in Trefftz polynomial spaces were proven. 
The derivation of $p$-version approximation estimates in higher dimensions is, to the best of our knowledge, still open. 
Here, we compare the Trefftz space to a full discontinuous polynomial space. 
Recall that, thanks to the Trefftz property of the test functions, we were able to cancel the volume integral in the weak formulation \eqref{eq:Trefftzweakform}. 
This does not hold for the full polynomial space. 
Therefore, we need to add the volume term to the left-hand side of the formulation \eqref{eq:TDG}, giving the new left-hand side:
\begin{equation*}
    \begin{aligned}
        & \mathcal{\tilde A}(v_{hp},\bm\sigma_{hp};w,\bm\tau ) := \\
        & \quad -\sum_{K\in\mathcal T_h}\int_K v_{hp} \left(\nabla\cdot \bm\tau + c^{-2} \dt w \right) + \bm\sigma_{hp}\cdot\left(\dt{\bm\tau} + \nabla w \right)\ dV \\
        & \quad +\mathcal{A}(v_{hp},\bm\sigma_{hp};w,\bm\tau ).
    \end{aligned}
\end{equation*}
We now need to solve $\mathcal{\tilde A}(v_{hp},\bm\sigma_{hp};w,\bm\tau )=\ell(w,\bm\tau),\ \forall (w,\bm\tau)\in\mathbb P^p(\mathcal T_h)^{n+1}$ for $(v_{hp},\bm\sigma_{hp})\in \mathbb P^p(\mathcal T_h)^{n+1}$.
Notice that, as opposed to the Trefftz-DG method \eqref{eq:TDG} the (full) DG method requires the computation of integrals also in space-time volumes.
\begin{figure*}[!ht]
    \centering
    \begin{subfigure}[t]{0.5\textwidth}
        \centering
        \includegraphics[width=\textwidth]{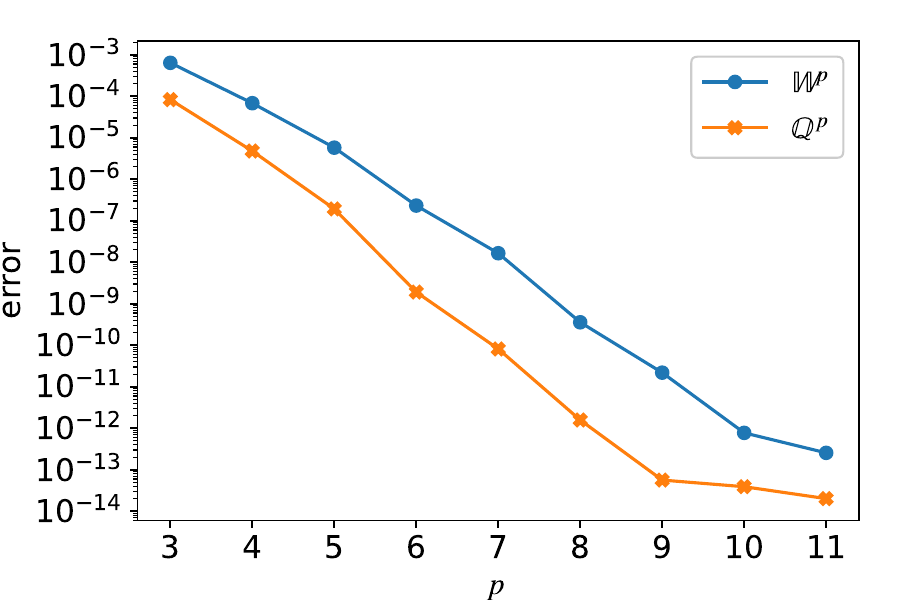}
    \end{subfigure}%
    \begin{subfigure}[t]{0.5\textwidth}
        \centering
        \includegraphics[width=\textwidth]{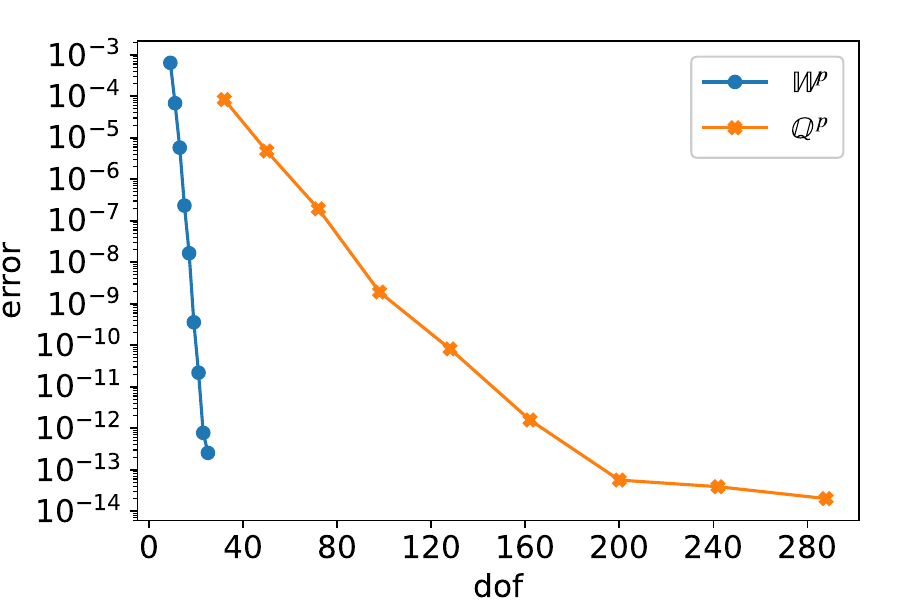}
    \end{subfigure}
    \caption{Comparison between Trefftz functions $\mathbb{W}^p$ and full polynomial space $\mathbb{Q}^p$, in terms of order (left) and local dofs (right).}
    \label{fig:pvt}
\end{figure*}
\par For the numerical results we have taken the unit square domain, uniformly meshed into space-time squares, and used the Trefftz space $\mathbb W^p$ or the space $\mathbb Q^p$ of polynomials with maximal order $p$ in each component.
%On squares it is only natural to construct the polynomial basis from a tensor product.
%Thus, we are in fact working with polynomials with maximal order $p$ in each component, which we denote by $\mathbb Q^p$. %i.e. $(v_{hp},\bm\sigma_{hp})\in (\mathbb P^p(\mathbb R) \otimes\mathbb P^p(\mathbb R) ) ^{1+1}$.
In the latter case the number of degrees of freedom per element is hence given by $2(p+1)^2$. 
\par The results are shown in \Cref{fig:pvt}. 
Both choices exhibit similar, exponential, convergence speed in terms of polynomial degree, 
although the Trefftz space is only a subset of the polynomials of maximal degree equal to $p$. 
%This is quite surprising, as the full polynomial space has a maximal degree of $p$ in each component, whereas the Trefftz space only uses polynomials of maximal polynomial degree equal to $p$. 
The benefits of the Trefftz space becomes clear when comparing errors versus number of degrees of freedom per element, as seen on the right in \Cref{fig:pvt}. 
%The Trefftz space uses much fewer degrees of freedom. 

\subsection{Comparing space-time meshing strategies}
\begin{figure*}[!ht]
    \centering
    \begin{subfigure}[t]{0.45\textwidth}
        \centering
        \includegraphics[width=\textwidth]{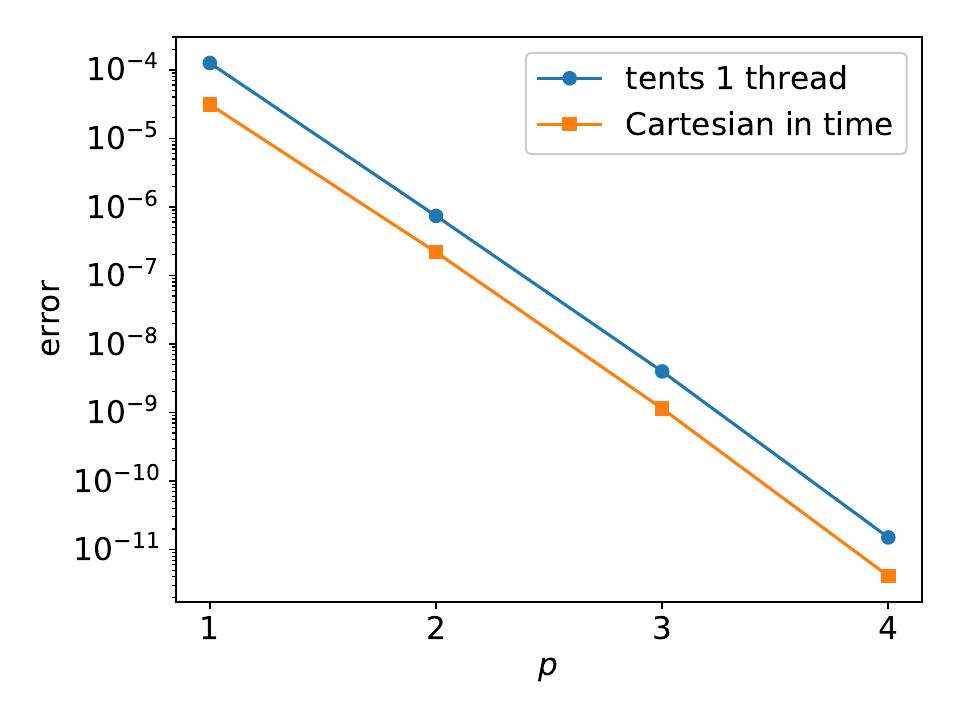}
    \end{subfigure}
    \begin{subfigure}[t]{0.45\textwidth}
        \centering
        \includegraphics[width=\textwidth]{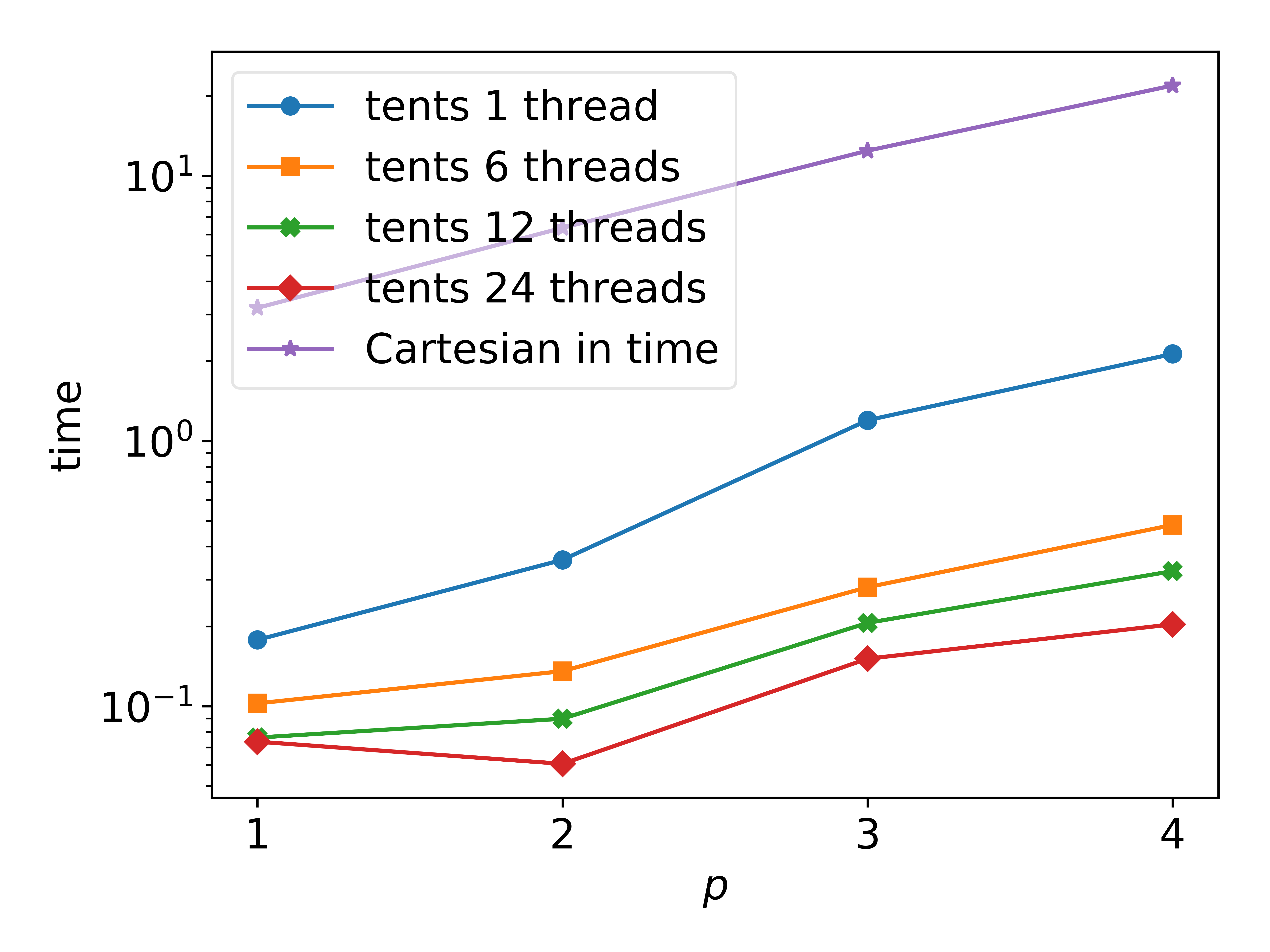}
        %\caption{}
    \end{subfigure}%
    \\
    \begin{subfigure}[t]{0.45\textwidth}
        \centering
        \includegraphics[width=\textwidth]{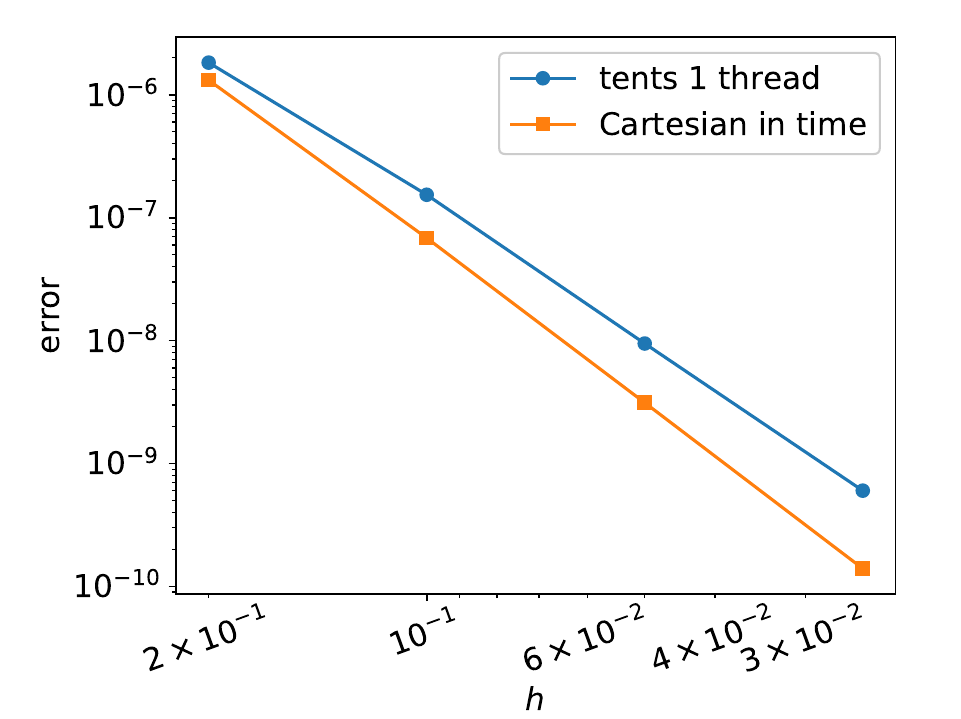}
    \end{subfigure}
    \begin{subfigure}[t]{0.45\textwidth}
        \centering
        \includegraphics[width=\textwidth]{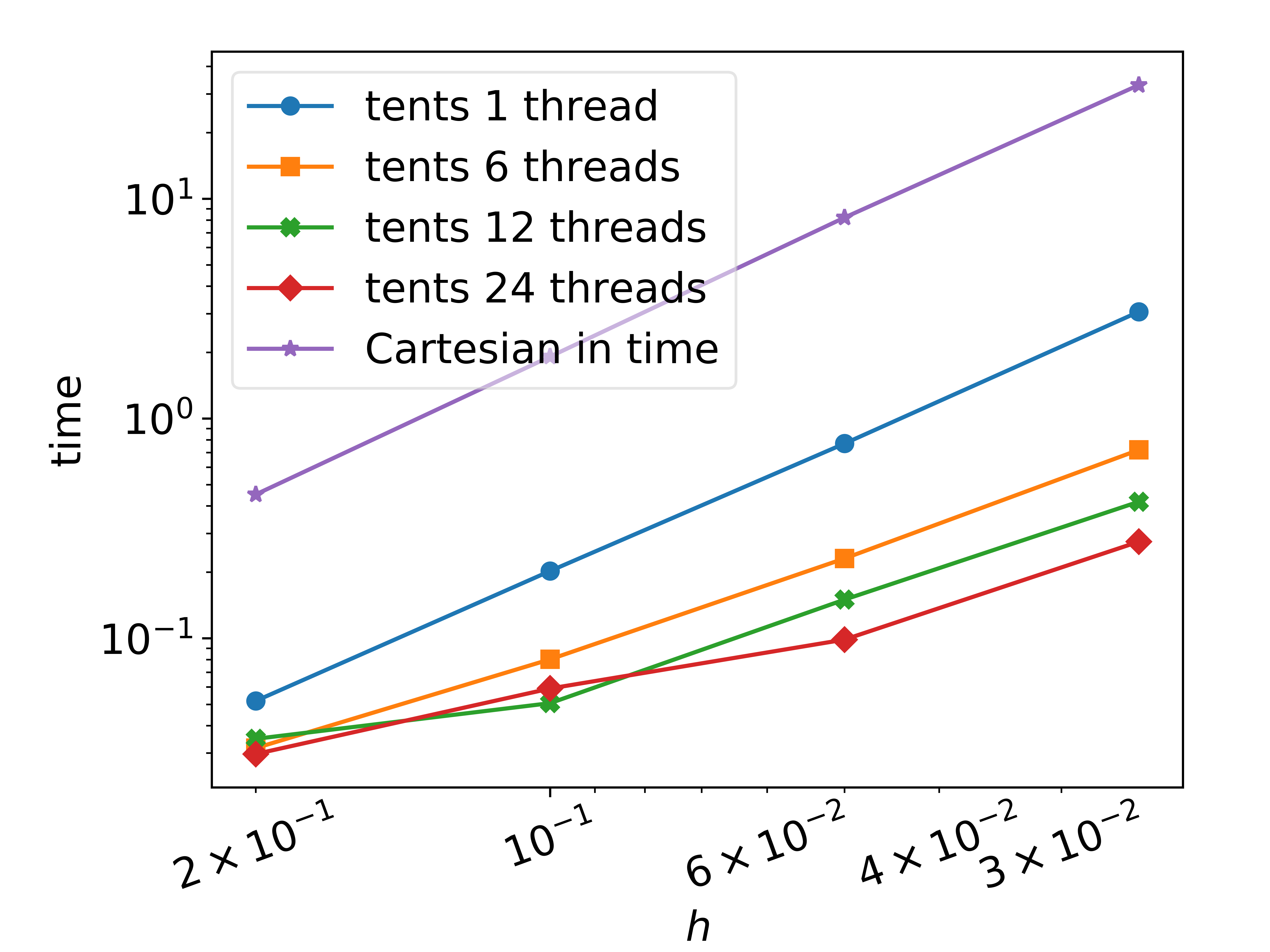}
        %\caption{}
    \end{subfigure}%
    \caption{Comparison of the Trefftz-DG method on Cartesian (in time) meshes and tent pitched meshes.} 
    \label{fig:tentvslabs}
\end{figure*}

In \Cref{sec:evoltent}, we have seen how to advance the solution element wise on a tent pitched mesh.
We now compare this approach to solving the full system on a Cartesian (in time) space-time slab. 
To solve the full system we use a block Jacobi solver.
When comparing the timing of the two methods, we consider 4 different cases for the tent pitching approach, first solving the tents sequentially, and then solving them in parallel on 6, 12, and 24 threads.
For this comparison, we choose a quasi-uniform mesh of the unit square in space and the final time equal to the mesh size, i.e. one CFL-conforming time step on the Cartesian mesh.
For the $p$-version comparison in \Cref{fig:tentvslabs} on the top we fix $h=0.04$, and in turn, for the $h$-version comparison on the bottom we fix $p=3$.
\par The results, in \Cref{fig:tentvslabs} on the left, show that the error between the two mesh types differs only slightly.
On the right in \Cref{fig:tentvslabs}, we compare the runtime of solving the full system on a Cartesian mesh, with solving the tents sequentially (on 1 thread), and solving them in parallel. 
Sequential tent pitching is about one magnitude faster.
Moreover, we can investigate the effects of parallelising the computations.
Using more threads only gives an advantage for small enough mesh sizes, as we are only able to solve independent tents in parallel. 

\subsection{Choice of spatial basis functions}\label{sec:compbasis}
As we have seen in \Cref{sec:Trefftzdisc}, the recursion formula \eqref{eq:Trefftz2ndorderrecursion} for the derivation of the Trefftz basis functions, can be initialized with an arbitrary choice of polynomial basis functions in space variable only. 
\begin{figure*}[!ht]
    \centering
    \begin{subfigure}[t]{0.5\textwidth}
        \centering
        \includegraphics[width=\textwidth]{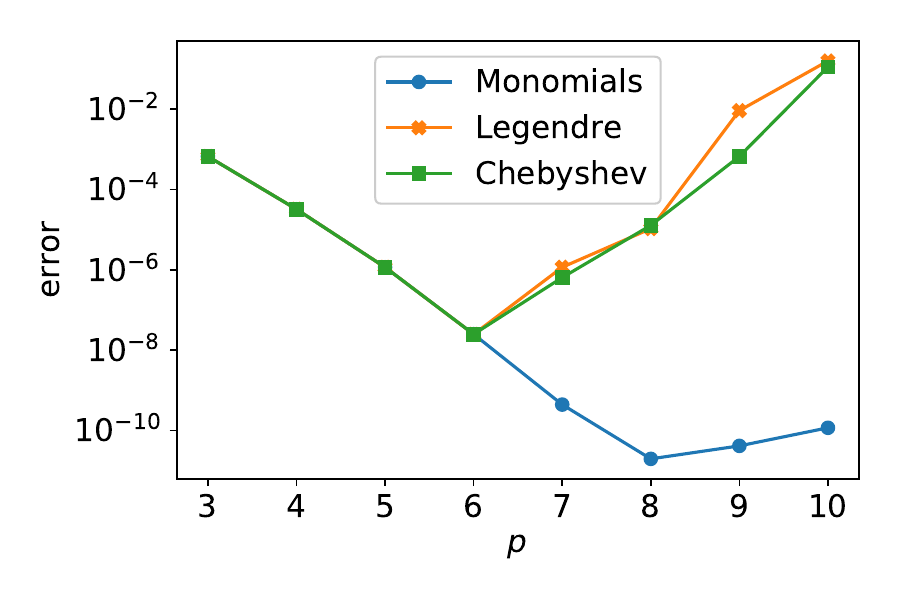}
        %\caption{L2 Error}
    \end{subfigure}%
    \begin{subfigure}[t]{0.5\textwidth}
        \centering
        \includegraphics[width=\textwidth]{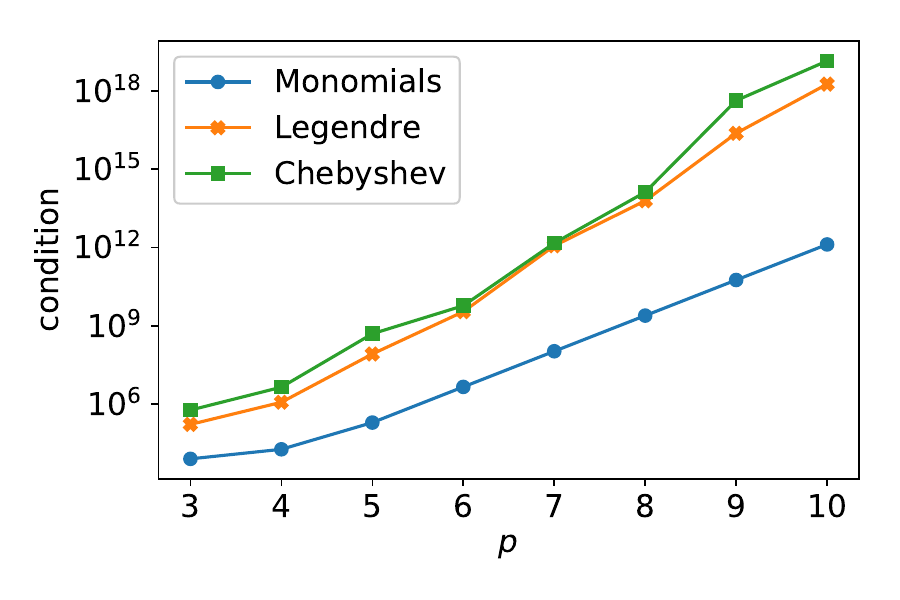}
        %\caption{Conditioning}
    \end{subfigure}
    \caption{Different types of initial polynomial basis functions. Comparison of the error (left) and the conditioning of the system global matrix (right).}
    \label{fig:basiscomp}
\end{figure*}
\par In the following, we compare three different choices for the initial polynomial basis functions: monomials, Legendre, and Chebychev polynomials. 
\ps{In all cases, the basis functions %Always, all polynomials
  are shifted to the center of each element and scaled, as described in \Cref{rem:scalepol}.}
We compare them in 1+1 dimensions, on the space-time unit square. 
The mesh considered is the tent pitched mesh shown in \Cref{fig:tentsol}. 
The problem is solved globally using formulation \eqref{eq:TDG}.
The results in \Cref{fig:basiscomp} show that all choices behave the same for low degrees. 
However, for higher degrees, Legendre and Chebychev polynomials fail to approximate the solution, 
%Legendre and Chebyshev show equally 
due to the bad conditioning of the system matrix, compared to the monomials. 
The good properties of the two sets of basis functions do not carry over when developed in the recursion.

\subsection{Tent pitching in 2 and 3 space dimensions}\label{sec:tp23}
\begin{figure*}[ht!]
    \centering
    \begin{subfigure}[t]{0.5\textwidth}
        \centering
        \includegraphics[width=\textwidth]{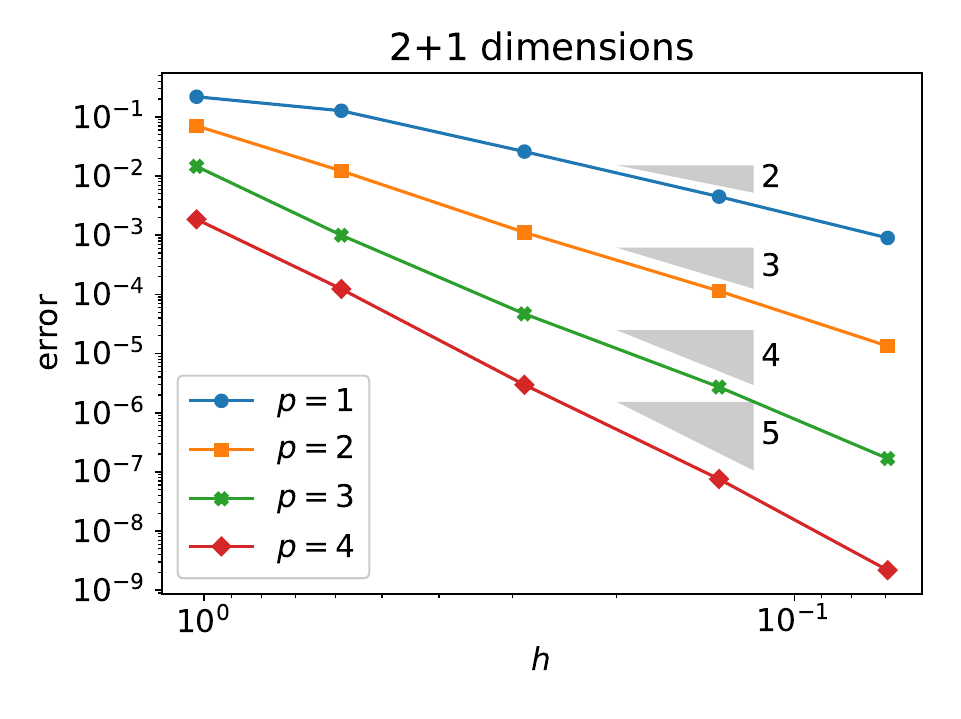}
        %\caption{}
    \end{subfigure}%
    \begin{subfigure}[t]{0.5\textwidth}
        \centering
        \includegraphics[width=\textwidth]{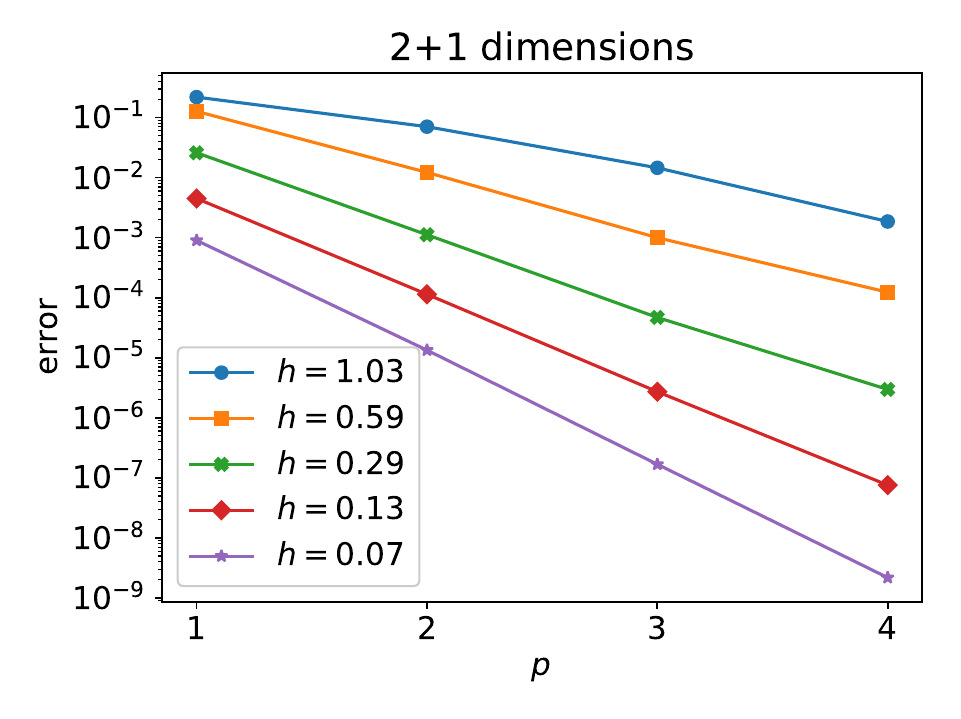}
    \end{subfigure}
    \\
    \centering
    \begin{subfigure}[t]{0.5\textwidth}
        \centering
        \includegraphics[width=\textwidth]{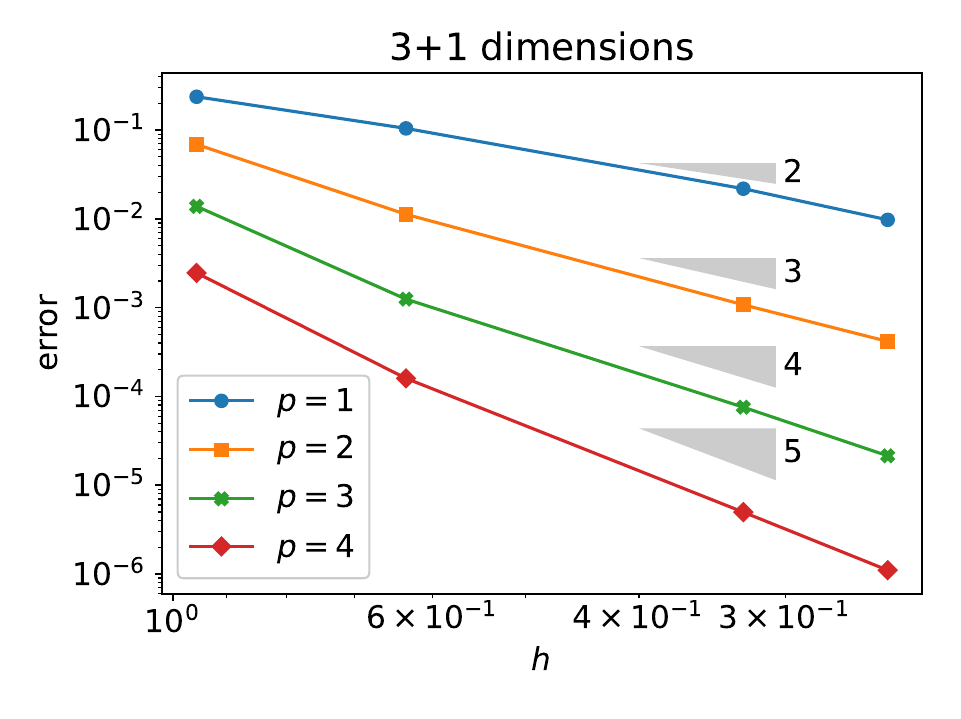}
    \end{subfigure}%
    \begin{subfigure}[t]{0.5\textwidth}
        \centering
        \includegraphics[width=\textwidth]{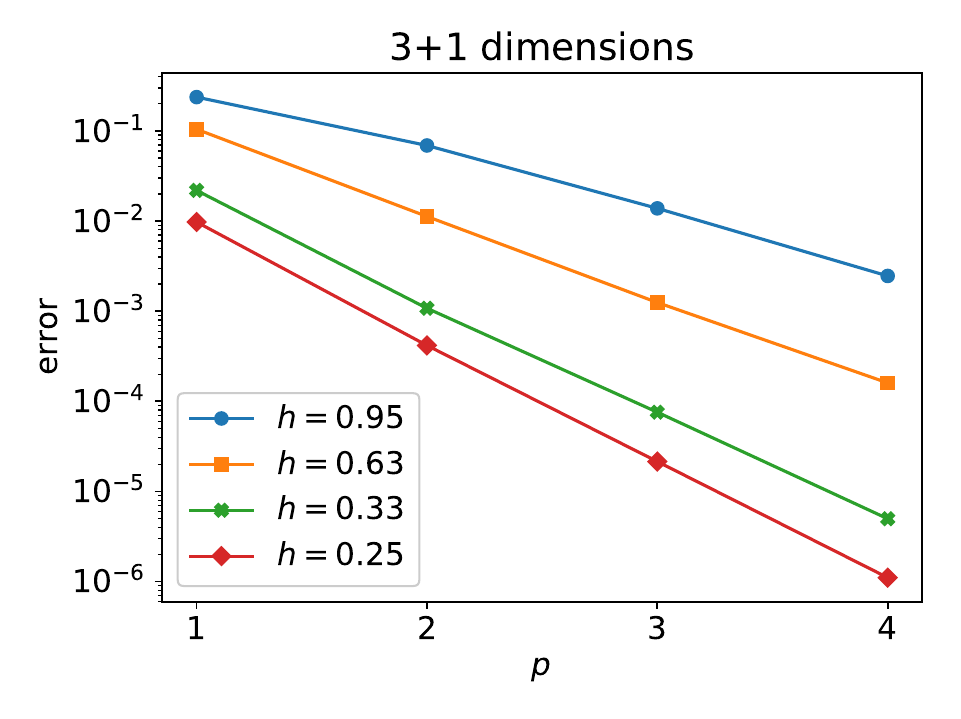}
    \end{subfigure}
    \caption{Tent pitching in 2+1 and 3+1 dimensions on the top and bottom, respectively. Convergence comparison with respect to the maximum mesh size on the left, and with respect to the Trefftz polynomial degree on the right. }\label{fig:tentpitching}
\end{figure*}

As discussed in \Cref{sec:evoltent}, we solve elementwise, and in parallel.
%Initial and Dirichlet boundary conditions are chosen, such that the standing wave is the solution. 
For this example, we choose as a spatial domain $\Omega$ the unit square and the unit cube. 
The initial quasi-uniform spatial mesh consists of triangles or tetrahedrons of maximal size $h$. 
We then use tent pitching in $2+1$ and $3+1$ dimensions, until the algorithm stops at time $T=1$, where we compute the error. 
The results of this are shown in \Cref{fig:tentpitching}. 
%In \cite{Moiola2017} it was shown that for a regular enough solution $(v,\bm\sigma)$ of the wave equation, $\|(v-v_{hp},\bm\sigma-\bm\sigma_{hp}\|\leq\mathcal O(h^{p+1/2})$.
\par In \Cref{fig:tentpitching} on the left, we plot the error in terms of $h$ for different values of polynomial degree $p$ of $\mathbb W ^p(Q)$. 
As typical for DG methods, in the case of a regular enough solution, we 
%can expect a convergence rate of $\mathcal O(h^{p+1/2})$, which 
%is confirmed by the numerical experiments. 
%is surpassed by the numerical examples, exhibiting a rate of $\mathcal O(h^{p+1})$, as typical for DG methods.
observe superconvergence, with the rate $\mathcal O(h^{p+1})$, outperforming the expected $\mathcal O(h^{p+1/2})$, see \cite[Thm. 6.19]{Moiola2017}.
%however we observe for larger mesh sizes the convergence speed cannot be observed immediately, especially in higher dimensions. \par
We also consider convergence in terms of degree $p$ of the Trefftz space $\mathbb W_p(\mathcal T_h)$, and report the results in \Cref{fig:tentpitching}, right plots.
For our analytic solution, we can observe exponential convergence. 

\subsection{Dissipation of energy}\label{sec:energy}
For smooth enough functions $(w,\bm\tau)$ the energy at a fixed time $\hat t$ is given by 
$$E(w,\bm\tau)=\frac12 \int_{\Omega\times\hat t} \left( c^{-2}w^2+|\bm\tau|^2 \right)\ dS.$$
In \cite{Moiola2017}, the method \eqref{eq:TDG} was shown to be dissipative, which we can also observe in numerical examples.
We test on a model problem with analytical solution
\begin{equation*}
    U(x,t)=\sin(\pi x)\sin(\pi t),
\end{equation*}
on the domain $[0,1]\times[0,T]$. 
We solve using the tent pitching algorithm. 
%An example solution up until $T=4$ is shown in \Cref{fig:tentslabsol}.
%\begin{figure*}[ht]
%\centering
%\includegraphics[width=0.4\textwidth]{../results/tentslabsol.jpg}
%\caption{Plot of the numerical solution on a tent pitched mesh in space-time, up to time $T=4$.}\label{fig:energy}  \label{fig:tentslabsol}
%\end{figure*}
\begin{figure*}[ht]
    \centering
    \begin{subfigure}[t]{0.5\textwidth}
        \centering
    \includegraphics[width=\textwidth]{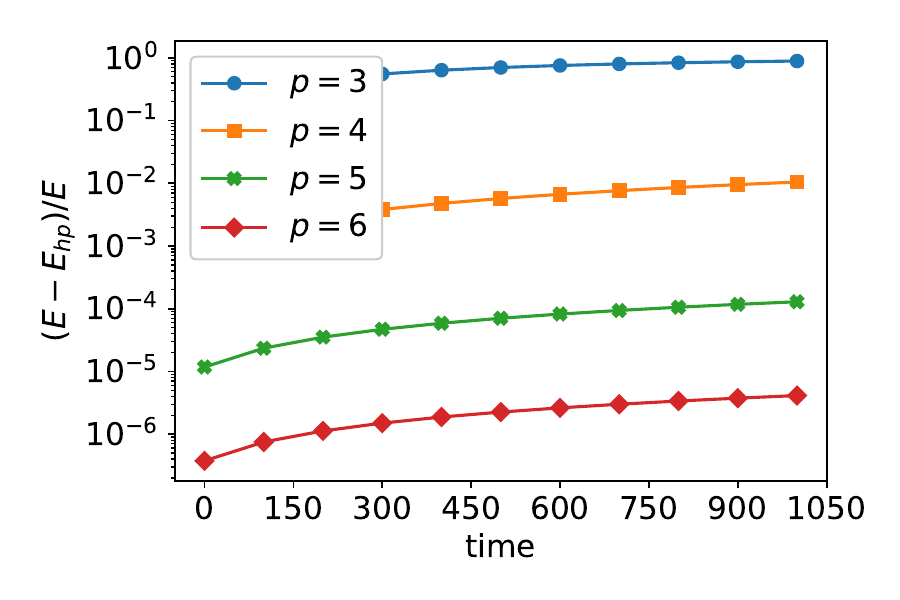}
        %\caption{}
    \end{subfigure}%
    \begin{subfigure}[t]{0.5\textwidth}
        \centering
    \includegraphics[width=\textwidth]{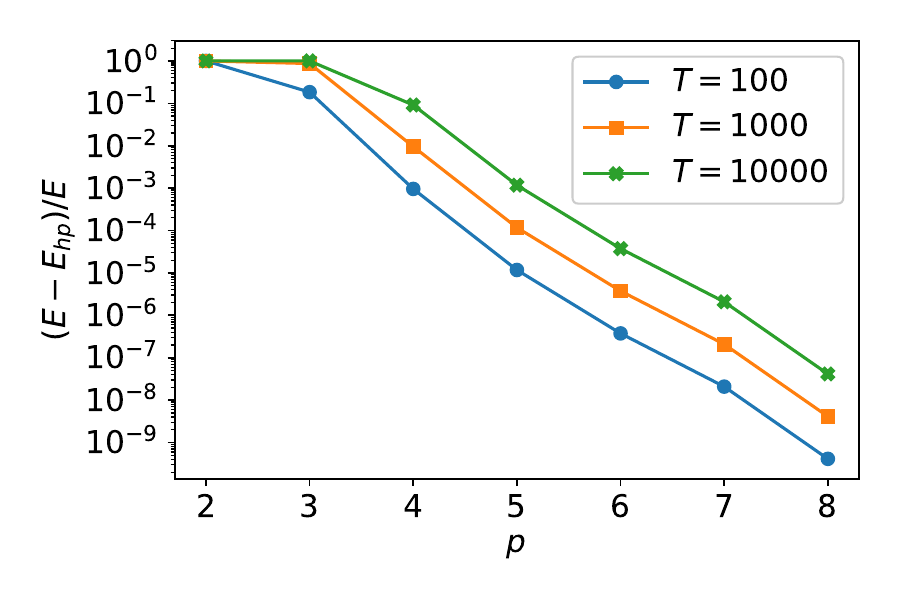}
    \end{subfigure}
    \caption{Error in energy over time and for different order of Trefftz polynomials.} \label{fig:energy}
\end{figure*}
The space mesh considered is a uniform partition of the interval $[0,1]$ into 5 elements.
We measure the relative error in the energy given by
\begin{equation*}
    \frac{E(\dt{ U},-\nabla U)-E(u_{hp},\bm\sigma_{hp})}{E(\dt{ U},-\nabla U)}.
\end{equation*}
%The results, shown in \Cref{fig:energy}, show the decay in energy, as expected. 
%In the figure we can see the results for Trefftz polynomials of order 5 and 6. 
%For higher order the decay in energy is slower. 
In \Cref{fig:energy} on the left, we can see that the error in energy
increases in time. %due to the dissipation of the method.
\ps{As the energy of the analytical solution is constant, and
 we are plotting the error without absolute value, we deduce that
the energy of the numerical solution is decreasing in time.}
The results actually suggest that the energy of the numerical solution decreases linearly in time.
In \Cref{fig:energy} on the right, we compare the error in the energy at three different times $T=100,1000,10000$, plotting it against the degree of Trefftz polynomials in a range from 2 to 8.
We observe exponential convergence for increasing order.
Furthermore, greater times $T$ seem to affect the error only by a multiplicative factor.

\subsection{Non-uniformly refined spatial meshes}\label{sec:lshape}
Now that we have verified the convergence of the Trefftz-DG method with tent pitching initialized on quasi-uniform spatial meshes, we test the advantages of the method on a non-uniformly refined spatial mesh.

\begin{table}[ht]
    \centering
    \begin{tabular}{|r|r|r|r|r|r|}
        \hline
        mesh   & $h_{\max}$ & total \#dofs& L2-error& dof-rate  & runtime [s]\\ \hline\hline
        \multirow{4}{*}{uniform} &  0.07  &  \num[round-precision=2,round-mode=figures, scientific-notation=true]{320112}  &  \num[round-precision=2,round-mode=figures, scientific-notation=true]{0.0178}  & - &  0.6234\\ \cline{2-6}
                             &  0.05  &  \num[round-precision=2,round-mode=figures, scientific-notation=true]{875696}  &  \num[round-precision=2,round-mode=figures, scientific-notation=true]{0.0117}  & 1.2638 &  1.5916\\ \cline{2-6}
                          &  0.03  &  \num[round-precision=2,round-mode=figures, scientific-notation=true]{3894752}  &  \num[round-precision=2,round-mode=figures, scientific-notation=true]{0.0067}  & 1.1022 &  7.9157\\ \cline{2-6}
 &  0.01  &  \num[round-precision=2,round-mode=figures, scientific-notation=true]{102320048}  &  \num[round-precision=2,round-mode=figures, scientific-notation=true]{0.0021}  & 1.0662 &  255.3233\\ \hline
    \multirow{4}{*}{non-uniform} &  0.12  &  \num[round-precision=2,round-mode=figures, scientific-notation=true]{1138256}  &  \num[round-precision=2,round-mode=figures, scientific-notation=true]{0.0218}  & - &  3.276\\ \cline{2-6}
                        &  0.10  &  \num[round-precision=2,round-mode=figures, scientific-notation=true]{2003536}  &  \num[round-precision=2,round-mode=figures, scientific-notation=true]{0.0083}  & 5.1308 &  4.6809\\ \cline{2-6}
                         &  0.08  &  \num[round-precision=2,round-mode=figures, scientific-notation=true]{3832624}  &  \num[round-precision=2,round-mode=figures, scientific-notation=true]{0.0030}  & 4.7041 &  8.0069\\ \cline{2-6}
 &  0.06  &  \num[round-precision=2,round-mode=figures, scientific-notation=true]{9777536}  &  \num[round-precision=2,round-mode=figures, scientific-notation=true]{0.0008}  & 4.2104 &  23.4588\\ \hline
    \end{tabular}
    \caption{Convergence rates and run time comparison for a singular solution on the L-shapes domain, comparing uniform meshing and meshes refined towards the singularity.}\label{tab:lshape}
\end{table}

In this test, the refinement is applied to resolve a singular solution at the reentrant corner of an L-shaped domain, given by $\Omega=[-1,1]^2\setminus([0,1]\times[-1,0])$. 
The mesh refinement strategy used takes the diameter of a spatial mesh elements $K$ as 
$$ h_K=h_{\max} r^{1-\mu},$$
where $r$ is the distance of $K$ to the reentrant corner, fixing a minimal mesh size of $h_{\min}=h_{\max}^{1/\mu}$.
Motivated by the theoretical results in \cite{apel}, we choose $\mu=\frac13$.

\par We consider a model problem with solution given, in polar coordinates, by
\begin{equation}\label{eq:vertgauss}
    U(r,\phi,t)=\cos(a t)\sin(\nu\phi)J_{\nu}(a r),
\end{equation}
where $J_{\nu}$ denotes the Bessel function of the first kind.
We consider $\nu=2/3$, so that $\nabla U$ is singular at the origin.
We solve up to time $T=1$ for $a=10$.
To avoid numerically integrate the singularity, we use the method to reconstruct the second order solution $U_{hp}\in \mathbb U^p(Q)$, introduced in \Cref{sec:recovery}, to measure the error given by $\|U(\cdot,T)-U_{hp}(\cdot,T)\|_{L^2(\Omega)}$.

\begin{figure*}[ht]
    \centering
    \begin{subfigure}[t]{0.35\textwidth}
        \centering
        \includegraphics[width=\textwidth]{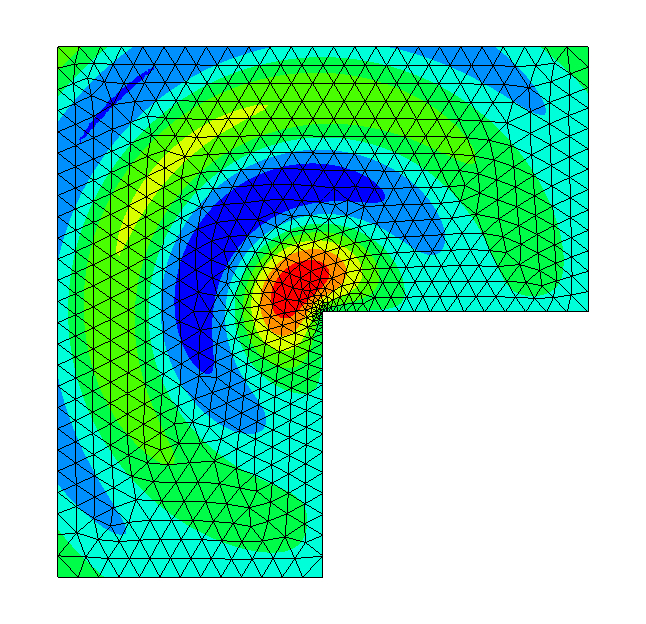}
    \end{subfigure}
    \begin{subfigure}[t]{0.5\textwidth}
        \centering
        \includegraphics[width=\textwidth]{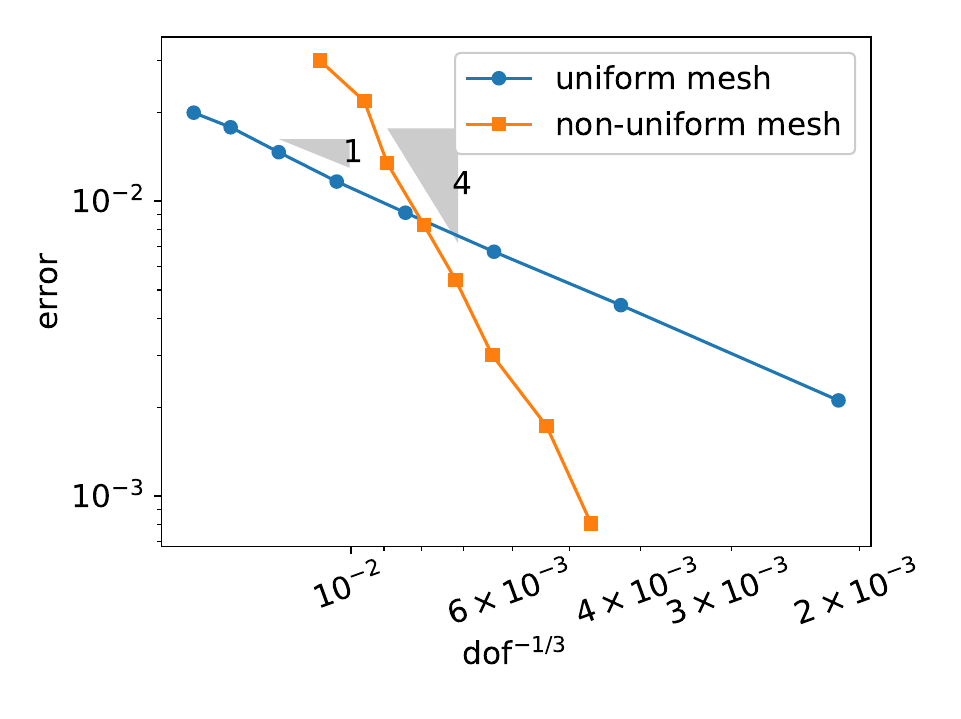}
        %\caption{}
    \end{subfigure}%
    \begin{subfigure}[t]{0.5\textwidth}
        \centering
        \includegraphics[width=\textwidth]{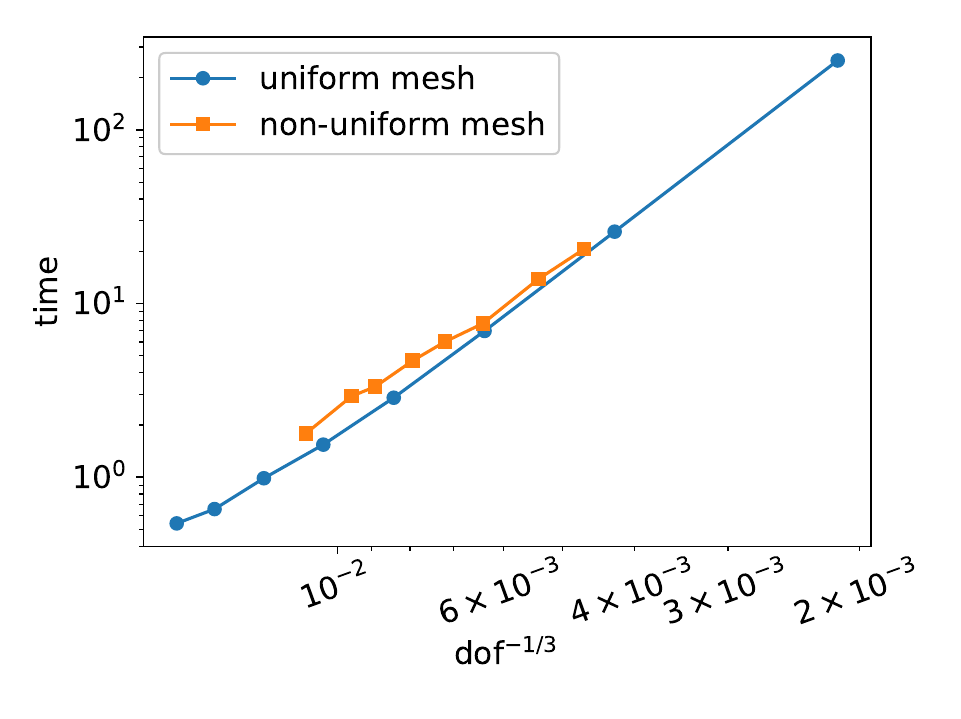}
    \end{subfigure}%
    \caption{The convergence rates on uniform and non-uniform meshes and timings (bottom) with the singular initial condition (top).}
    \label{fig:lshape}
\end{figure*}

\par The comparison between results obtained with uniform and non-uniform mesh refinement are shown in \Cref{fig:lshape} for Trefftz functions of degree $p=3$. 
We compare the two different meshing strategies by plotting them against (global dof)$^{-1/3}$. 
For the uniformly refined meshes, the convergence rate is bounded by the smoothness of the solution $U\in H^{5/3-\varepsilon}(Q)$, for $\varepsilon>0$. We observe a convergence rate of $\mathcal O (h^{1})$. 
%This is observed in the results shown in \Cref{fig:lshape}.
Using the non-uniformly refined meshes, we are able to recover optimal convergence for the third order Trefftz polynomials, as seen in \Cref{fig:lshape}. 

\Cref{tab:lshape}, gives a closer look on some of the properties already visualized in \Cref{fig:lshape} and also shows the runtime (in seconds).
%The computations were performed on a Lenovo X230, Intel i7-3520M CPU @ 2.90Ghz. 
For the computations we used 24 threads.
In \Cref{fig:lshape} on the bottom right, we compare the run time with the degrees of freedom. 
We observe that the uniform and the non-uniform mesh take about the same time for comparable numbers of degrees of freedom. 
Thus, no significant locking, due to the spatial refinement, occurs. 

\subsection{Wave propagation in an heterogeneous material}\label{sec:mat}
In the following example we investigate the reflection of a wave at an interface of two different materials. This experimental setup was also performed in \cite{Kocher, Bangerth_adaptivegalerkin}.
We consider the space-time domain $Q=[0,2]^2\times(0,1]$, and problem \eqref{eq:TDG} with homogeneous Dirichlet boundary conditions.
%The domain is split into two parts with different material properties, given by the piecewise constant wavenumber
The wavespeed is the piecewise constant function given by
\begin{align*}
    c(x_1,x_2) =    
    \begin{cases}
        1 & x_1\leq1.2,\\    
        3 & x_1>1.2.
    \end{cases}
\end{align*}
\begin{figure*}[!hb]
    \centering
    \begin{subfigure}[t]{0.5\textwidth}
        \centering
        \includegraphics[width=\textwidth]{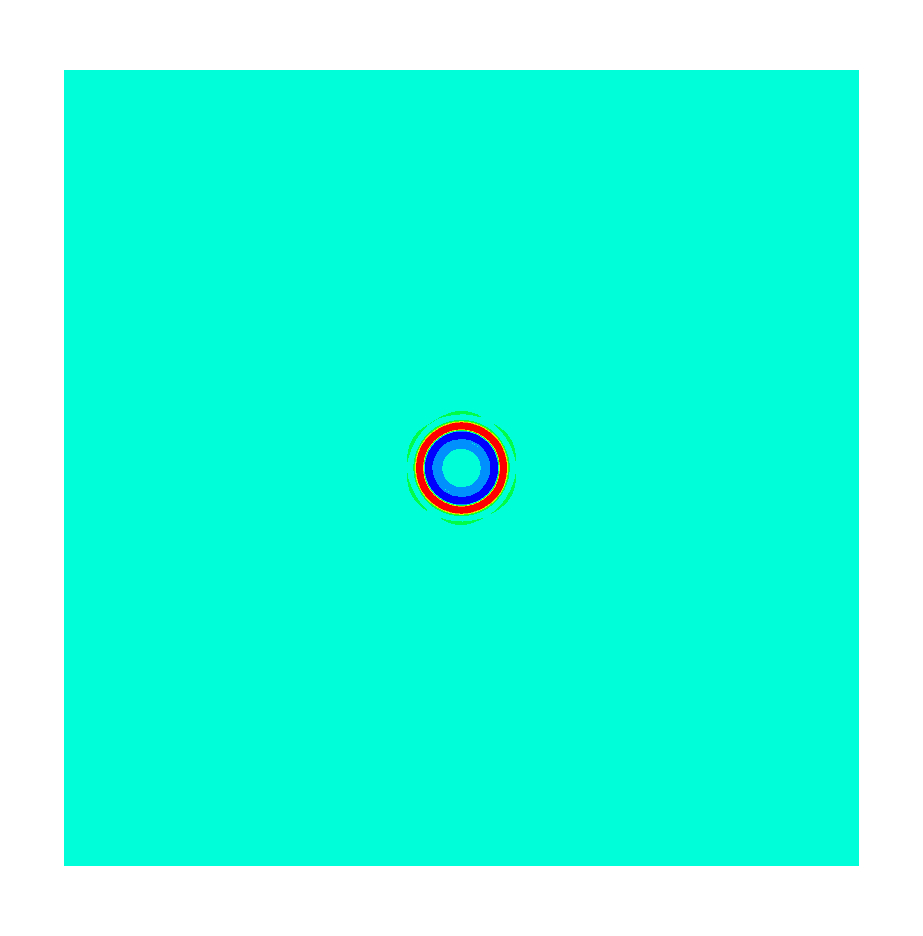}
        %\caption{}
    \end{subfigure}%
    \begin{subfigure}[t]{0.5\textwidth}
        \centering
        \includegraphics[width=\textwidth]{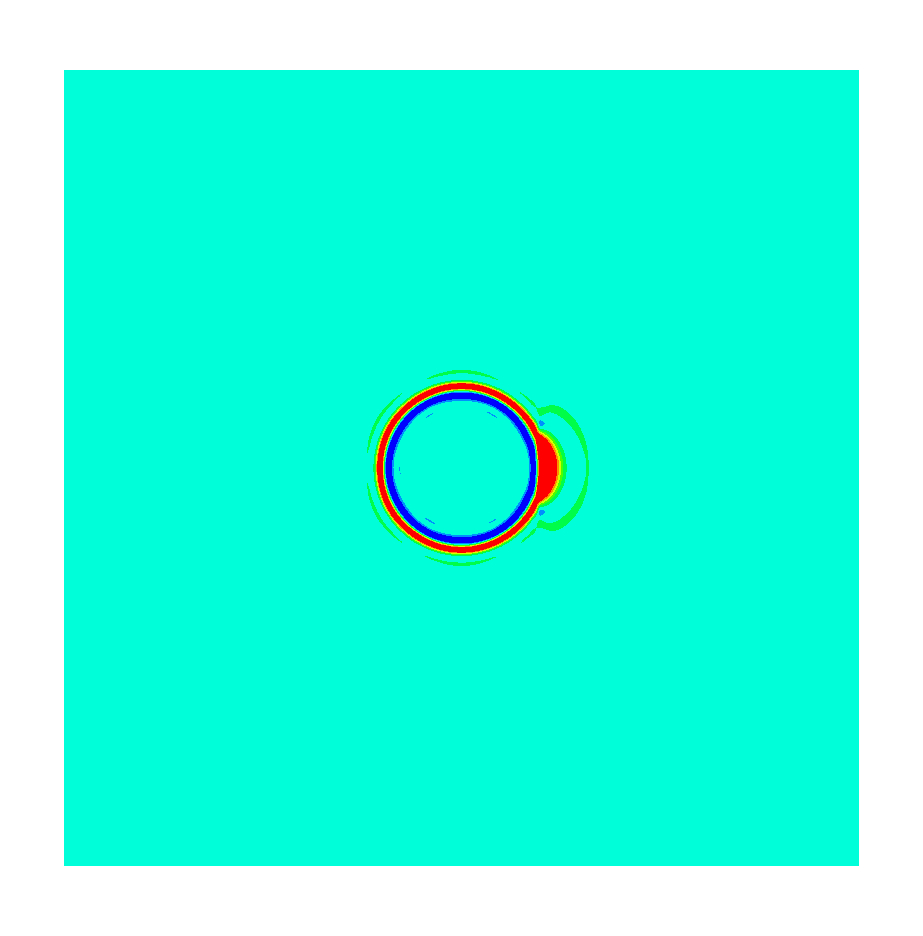}
    \end{subfigure}
    \\
    \centering
    \begin{subfigure}[t]{0.5\textwidth}
        \centering
        \includegraphics[width=\textwidth]{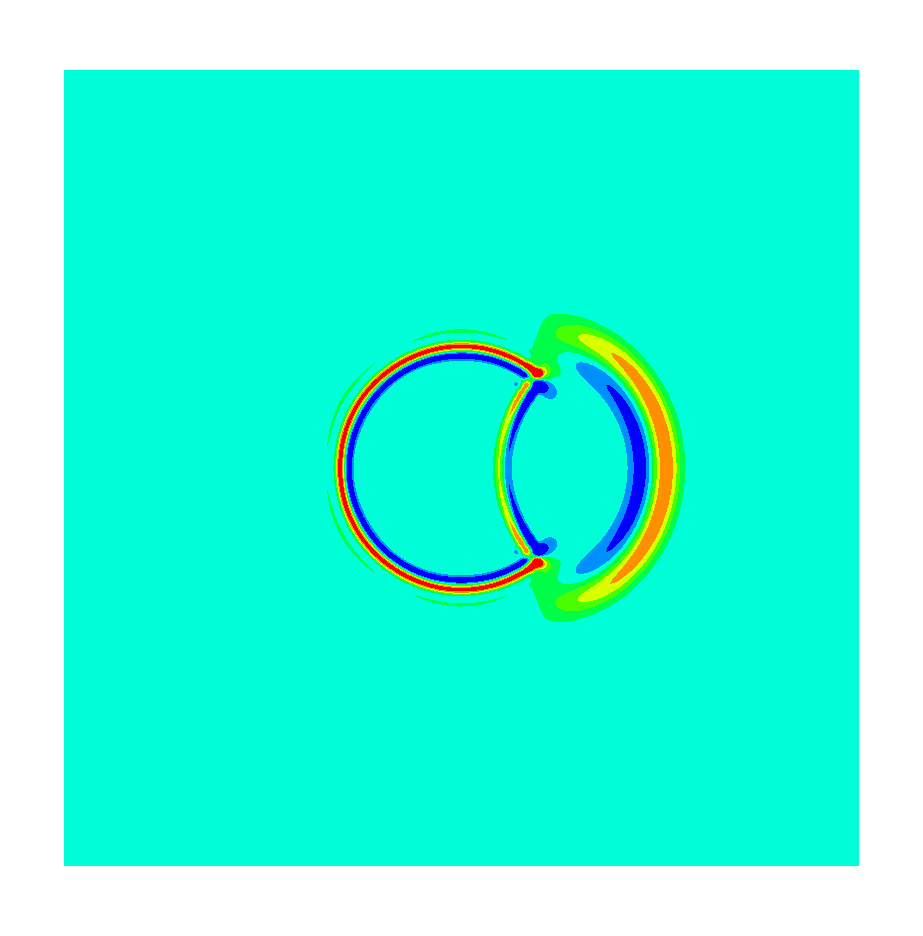}
        %\caption{}
    \end{subfigure}%
    \begin{subfigure}[t]{0.5\textwidth}
        \centering
        \includegraphics[width=\textwidth]{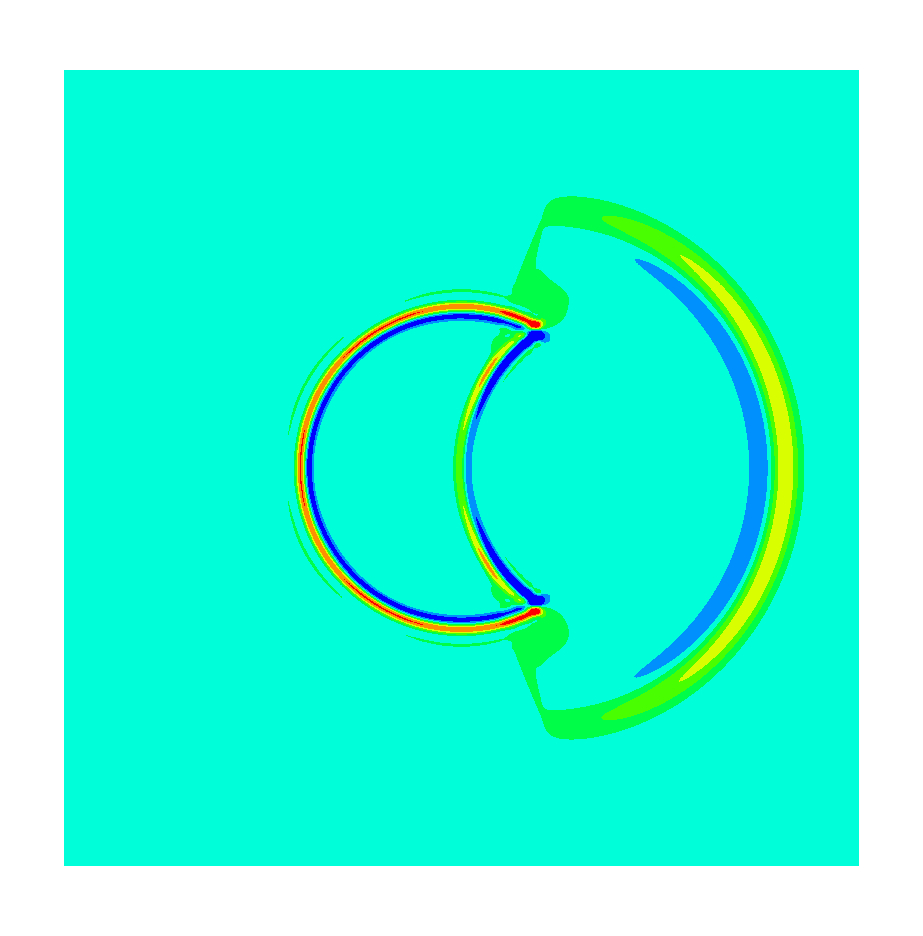}
    \end{subfigure}
    \caption{Wave traveling through inhomogeneous material, shown at times $T=0.1,0.2,0.3,0.4$.}\label{fig:solmat}
\end{figure*}\\
As initial condition, we take a Gaussian wave given by
\begin{align*}
    U_0(\bm x)=\exp{(-\|\bm x-\bm x_0\|^2/\delta^2)},\
    v_0(\bm x)=0,
\end{align*}
where we choose $\bm x_0=(1,1)$ and $\delta=0.01$.
The computations are performed with polynomial degree $p=4$.

\par Snapshots of the solution are shown in \Cref{fig:solmat}.
In the Snapshots, the right part of the domain has spatial mesh sizes up to $0.03$, whereas in the left part we choose as spatial mesh size of $0.01$, in order to better capture the steeper wavefront in the slower traveling material. 
First, we see that the initial condition unfolds in the left homogeneous part of the medium.
At $T=0.2$, the wave crosses over into the material with higher wave velocity. 
In the next snapshot we can see that the wave splits into a part traveling to the right with a higher velocity and shallow wavefront, and a part reflected at the interface traveling backwards to the left. 
Finally, at the time $T=0.4$, we can also observe the weaker Huygens wave, which traveled parallel to the interface, before traveling back towards the left.

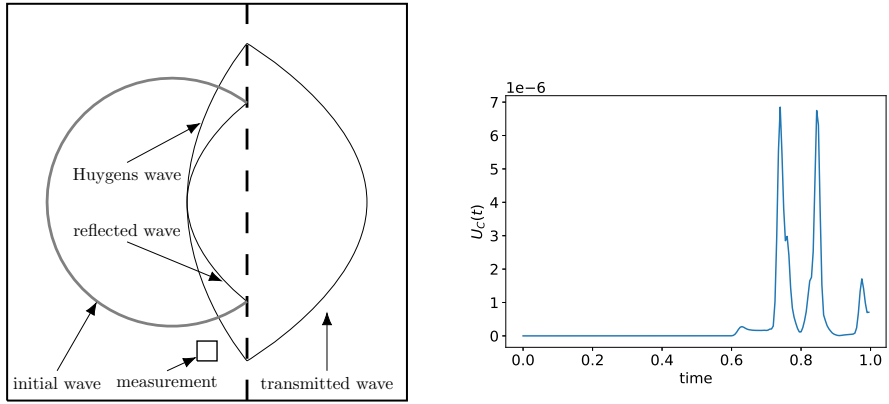
\begin{figure}[ht]\centering
    \centering
    %\begin{subfigure}[t]{0.49\textwidth}
        \centering
    \resizebox{0.5\linewidth}{!}{
    \begin{tikzpicture}[>=latex,scale=1, every node/.style={transform shape}]
        %draw thick boundary
        \draw[line width=0.1mm] (0,0) -- (2,0) -- (2,2) -- (0,2) -- (0,0);
        %%split domain
        \draw[dashed] (1.2,0) -- (1.2,2); 
        %%measure
        \draw[line width=0.07mm] (1-0.05,0.25-0.05) --(1+0.05,0.25-0.05) --(1+0.05,0.25+0.05) --(1-0.05,0.25+0.05) --(1-0.05,0.25-0.05);
        %waves
        \draw [line width=0.05mm] plot [smooth, tension=1] coordinates {(1.2,1.8) (1.8,1) (1.2,0.2)};
        \draw [line width=0.05mm] plot [smooth, tension=1] coordinates {(1.2,1.8) (0.9,1) (1.2,0.2)};
        \draw [line width=0.05mm] plot [smooth, tension=1] coordinates {(1.2,1.5) (0.9,1) (1.2,0.5)};
        \tkzDefPoint(1.2,1.5){A}\tkzDefPoint(0.2,1){B}\tkzDefPoint(1.2,0.5){C}
        \tkzCircumCenter(A,B,C)
        \tkzGetPoint{O}
        \tkzDrawArc(O,A)(C)
        %names
        \draw[line width=0.05mm,-{Latex[length=0.7mm]}] (1.6,0.15) node[anchor=north ,scale=0.25] {transmitted wave} -- (1.6,0.45) ;
        \draw[line width=0.05mm,-{Latex[length=0.7mm]}] (0.6,0.8) node[anchor=south ,scale=0.25] {reflected wave} -- (1.08,0.6) ;
        \draw[line width=0.05mm,-{Latex[length=0.7mm]}] (0.6,1.2) node[anchor=north ,scale=0.25] {Huygens wave}-- (0.97,1.4) ;
        \draw[line width=0.05mm,-{Latex[length=0.7mm]}] (0.25,0.15) node[anchor=north,scale=0.25] {initial wave}-- (0.45,0.5) ;
        \draw[line width=0.05mm,-{Latex[length=0.7mm]}] (0.8,0.15) node[anchor=north,scale=0.25] {measurement}-- (0.95,0.22) ;
        %\draw [line width=0.05mm] plot [smooth, tension=4] coordinates {(1.2,1.5) (0.2,1) (1.2,0.5)};
%\def\centerarc[#1](#2)(#3:#4:#5)% Syntax: [draw options] (center) (initial angle:final angle:radius)
    %{ \draw[#1] ($(#2)+({#5*cos(#3)},{#5*sin(#3)})$) arc (#3:#4:#5); }
%\centerarc[](-1,-1)(5:85:1)
       %\draw [domain=170:210] plot ({1+cos(\x)}, {1+sin(\x)});
    \end{tikzpicture}
    }
    \begin{subfigure}[t]{0.49\textwidth}
        \centering
        \includegraphics[width=\textwidth]{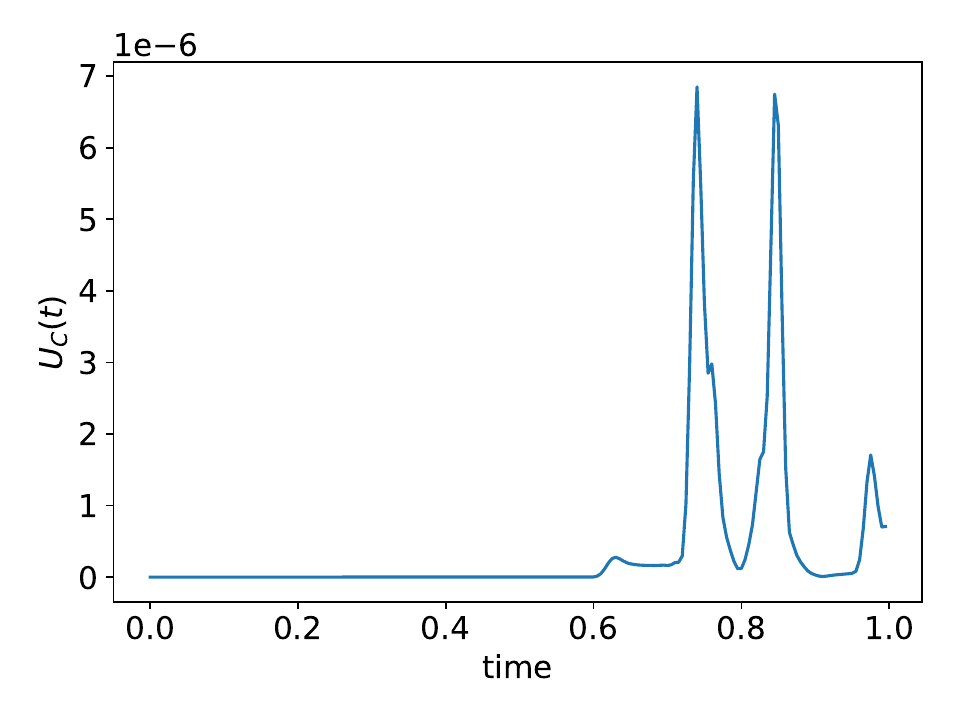}
    %\resizebox{\linewidth}{!}{ %}
    \end{subfigure}
%\end{subfigure}
    %\caption{A product mesh on the left and a tent pitched mesh on the right. Both are on top the same spatial mesh, and consist of vertical faces or faces below the characteristic speed.}
    \caption{Sketch of the expected wave pattern (left) and measured output quantity (right).}
    \label{fig:singu}
\end{figure}
\par In \Cref{fig:singu} on the left we present a sketch of the actions described above, also indicating a region where we measured the output
$$U_C(t)=\|U(\cdot,t)\|_{L^1(\Omega_C)}.$$
The domain of measurement was chosen $\Omega_C=[1-\varepsilon_C,1+\varepsilon_C]\times[0.25-\varepsilon_C,0.25+\varepsilon_C]$, with $\varepsilon_C=2^{-7}$. 
The measurement over time is presented in \Cref{fig:singu} on the right and shows that we are able to distinguish the three incoming waves. 
We can see the very weak Huygens wave arriving first, followed by the initial wave and the reflected one.

\section{Conclusion}
We have presented implementational aspects and numerical results for the Trefftz-DG method for the acoustic wave equation, originally presented in \cite{Moiola2017}. 
The implementation in NGSolve was used to solve on Cartesian (in time) meshes and tent pitched meshes in up to 3+1 dimensions, with varying mesh sizes, polynomial degrees, and wavenumber.
\par 
The $h$-convergence rates were shown to comply with the analytic results, showing superconvergence for analytic solutions and limited rates in the case of solutions with insufficient regularity. 
In the latter case, we were able to recover optimal convergence using non-uniform meshes.
For analytic solutions we observed exponential convergence rates in the polynomial degree $p$. 
%The $p$-convergence of the method is still an open problem. 
%, as well as the analysis on more general meshes (allowing slanted faces above the characteristic speed), are
\par Possible developments include the extension to the case of electromagnetic waves (Maxwell's equations).

%\section*{Acknowledgment}
\begin{acknowledgment*}
This work has been supported by the Austrian Science Fund, grants no. W1245 and F65.
\end{acknowledgment*}

\nocite{*}
%\small
\bibliographystyle{abbrv}
\bibliography{bibliografie}

\begin{thebibliography}{10}

\bibitem{apel}
T.~Apel.
\newblock {\em Anisotropic finite elements: local estimates and applications}.
\newblock Advances in Numerical Mathematics. B. G. Teubner, Stuttgart, 1999.

\bibitem{Bangerth_adaptivegalerkin}
W.~Bangerth, M.~Geiger, and R.~Rannacher.
\newblock Adaptive {G}alerkin finite element methods for the wave equation.
\newblock {\em Comput. Methods Appl. Math.}, 10(1):3--48, 2010.

\bibitem{Banjai}
L.~Banjai, E.~H. Georgoulis, and O.~Lijoka.
\newblock A {T}refftz polynomial space-time discontinuous {G}alerkin method for
  the second order wave equation.
\newblock {\em SIAM J. Numer. Anal.}, 55(1):63--86, 2017.

\bibitem{barucq}
H.~Barucq, H.~Calandra, J.~Diaz, and E.~Shishenina.
\newblock {Space-time Trefftz-DG approximation for elasto-acoustics}.
\newblock {\em {Appl. Anal.}}, 00:1 -- 16, Aug. 2018.

\bibitem{Drfler}
W.~D\"{o}rfler, S.~Findeisen, and C.~Wieners.
\newblock Space-time discontinuous {G}alerkin discretizations for linear
  first-order hyperbolic evolution systems.
\newblock {\em Comput. Methods Appl. Math.}, 16(3):409--428, 2016.

\bibitem{Egger2015TransparentBC}
H.~Egger, F.~Kretzschmar, S.~M. Schnepp, I.~Tsukerman, and T.~Weiland.
\newblock Transparent boundary conditions for a discontinuous {G}alerkin
  {T}refftz method.
\newblock {\em Appl. Math. Comput.}, 267:42--55, 2015.

\bibitem{EggerKretzMaxwell}
H.~Egger, F.~Kretzschmar, S.~M. Schnepp, and T.~Weiland.
\newblock A space-time discontinuous {G}alerkin {T}refftz method for time
  dependent {M}axwell's equations.
\newblock {\em SIAM J. Sci. Comput.}, 37(5):B689--B711, 2015.

\bibitem{tp2}
J.~Erickson, D.~Guoy, J.~Sullivan, and A.~Üngör.
\newblock {B}uilding space-time meshes over arbitrary spatial domains.
\newblock {\em Eng. Comput.}, 20:342--353, 2005.

\bibitem{FalkRichter}
R.~S. Falk and G.~R. Richter.
\newblock Explicit finite element methods for symmetric hyperbolic equations.
\newblock {\em SIAM J. Numer. Anal.}, 36(3):935--952, 1999.

\bibitem{MTPMaxwell}
J.~Gopalakrishnan, M.~Hochsteger, J.~Schöberl, and C.~Wintersteiger.
\newblock {A}n explicit mapped tent pitching scheme for {M}axwell equations,
  2019.
\newblock arXiv:1906.11029.

\bibitem{gopala}
J.~Gopalakrishnan, P.~Monk, and P.~Sep\'{u}lveda.
\newblock A tent pitching scheme motivated by {F}riedrichs theory.
\newblock {\em Comput. Math. Appl.}, 70(5):1114--1135, 2015.

\bibitem{MTP}
J.~Gopalakrishnan, J.~Sch\"{o}berl, and C.~Wintersteiger.
\newblock Mapped tent pitching schemes for hyperbolic systems.
\newblock {\em SIAM J. Sci. Comput.}, 39(6):B1043--B1063, 2017.

\bibitem{grote2}
M.~J. Grote, M.~Mehlin, and S.~A. Sauter.
\newblock Convergence analysis of energy conserving explicit local
  time-stepping methods for the wave equation.
\newblock {\em SIAM J. Numer. Anal.}, 56(2):994--1021, 2018.

\bibitem{grote1}
M.~J. Grote and T.~Mitkova.
\newblock High-order explicit local time-stepping methods for damped wave
  equations.
\newblock {\em J. Comput. Appl. Math.}, 239:270--289, 2013.

\bibitem{Hughes}
T.~J.~R. Hughes and G.~M. Hulbert.
\newblock Space-time finite element methods for elastodynamics: formulations
  and error estimates.
\newblock {\em Comput. Methods Appl. Mech. Engrg.}, 66(3):339--363, 1988.

\bibitem{Kocher}
U.~K\"{o}cher and M.~Bause.
\newblock Variational space-time methods for the wave equation.
\newblock {\em J. Sci. Comput.}, 61(2):424--453, 2014.

\bibitem{KretzMoiolaMaxwell}
F.~Kretzschmar, A.~Moiola, I.~Perugia, and S.~M. Schnepp.
\newblock {\it {A} priori} error analysis of space-time {T}refftz discontinuous
  {G}alerkin methods for wave problems.
\newblock {\em IMA J. Numer. Anal.}, 36(4):1599--1635, 2016.

\bibitem{LILIENTHAL}
M.~Lilienthal, S.~M. Schnepp, and T.~Weiland.
\newblock Non-dissipative space-time {$hp$}-discontinuous {G}alerkin method for
  the time-dependent {M}axwell equations.
\newblock {\em J. Comput. Phys.}, 275:589--607, 2014.

\bibitem{Lowrie}
R.~B. Lowrie, P.~L. Roe, and B.~van Leer.
\newblock Space-time methods for hyperbolic conservation laws.
\newblock In {\em Barriers and challenges in computational fluid dynamics
  ({H}ampton, {VA}, 1996)}, volume~6 of {\em ICASE/LaRC Interdiscip. Ser. Sci.
  Eng.}, pages 79--98. Kluwer Acad. Publ., Dordrecht, 1998.

\bibitem{Moiola2017}
A.~Moiola and I.~Perugia.
\newblock A space-time {T}refftz discontinuous {G}alerkin method for the
  acoustic wave equation in first-order formulation.
\newblock {\em Numer. Math.}, 138(2):389--435, 2018.

\bibitem{Pete-2005}
P.~Monk and G.~R. Richter.
\newblock A discontinuous {G}alerkin method for linear symmetric hyperbolic
  systems in inhomogeneous media.
\newblock {\em J. Sci. Comput.}, 22/23:443--477, 2005.

\bibitem{Richter}
G.~R. Richter.
\newblock An explicit finite element method for the wave equation.
\newblock {\em Appl. Numer. Math.}, 16(1-2):65--80, 1994.

\bibitem{joachim}
J.~Sch{\"o}berl.
\newblock {C++11 implementation of finite elements in NGSolve}.
\newblock {\em ASC Report 30/2014, Institute for Analysis and Scientific
  Computing, Vienna University of Technology}, 2014.

\bibitem{ngsolve}
J.~{Schöberl}.
\newblock {NGSolve finite element library}.
\newblock https://ngsolve.org.
\newblock Accessed: 2019-07-24.

\bibitem{Steinbach2019}
O.~Steinbach and M.~Zank.
\newblock {\em A stabilized space--time finite element method for the wave
  equation}, pages 341--370.
\newblock Springer International Publishing, Cham, 2019.

\bibitem{trefftz}
E.~Trefftz.
\newblock {Ein Gegenstück zum Ritzschen Verfahren}.
\newblock {\em Verhandl. 2er Internat. Kongress. Techn. Mechanik Zürich, 1926,
  12–17 Sept.}, pages 131--137, 1926.

\bibitem{tp1}
A.~\"{U}ng\"{o}r and A.~Sheffer.
\newblock Pitching tents in space-time: mesh generation for discontinuous
  {G}alerkin method.
\newblock {\em Internat. J. Found. Comput. Sci.}, 13(2):201--221, 2002.
\newblock Volume and surface triangulations.

\bibitem{WinterMTP}
C.~Wintersteiger.
\newblock Mapped tent pitching method for hyperbolic conservation laws.
\newblock {\em {D}iplomarbeit}, 2015.

\bibitem{10.1007/978-3-642-59721-3_48}
L.~Yin, A.~Acharya, N.~Sobh, R.~B. Haber, and D.~A. Tortorelli.
\newblock A space-time discontinuous {G}alerkin method for elastodynamic
  analysis.
\newblock In {\em Discontinuous {G}alerkin methods ({N}ewport, {RI}, 1999)},
  volume~11 of {\em Lect. Notes Comput. Sci. Eng.}, pages 459--464. Springer,
  Berlin, 2000.

\end{thebibliography}

\end{document}